# Lattice Structures for Attractors I

WILLIAM D. KALIES

Florida Atlantic University
777 Glades Road
Boca Raton, FL 33431, USA

KONSTANTIN MISCHAIKOW

Rutgers University
110 Frelinghusen Road
Piscataway, NJ 08854, USA

ROBERT C.A.M. VANDERVORST

VU University
De Boelelaan 1081a
1081 HV, Amsterdam, The Netherlands

ABSTRACT. We describe the basic lattice structures of attractors and repellers in dynamical systems. The structure of distributive lattices allows for an algebraic treatment of gradient-like dynamics in general dynamical systems, both invertible and noninvertible. We separate those properties which rely solely on algebraic structures from those that require some topological arguments, in order to lay a foundation for the development of algorithms to manipulate these structures computationally.

## 1. Introduction

As is made clear by Conley [**6**], attractors are central to our theoretical understanding of global nonlinear dynamics in that they form the basis for robust decompositions of gradient-like structures. They also play a singularly important role in applications in that they often represent the dynamics that is observed. Thus it is not surprising that attractors appear as standard topics in nonlinear dynamics [**15**], and that their structure has been studied in a wide variety of settings. We make no attempt to provide even a cursory list of references on the subject, but we do remark that the theory developed in [**1**] applies in a general setting in which the dynamics is generated by relations.

Our motivation for returning to this subject arises from computations and applications. With regard to applications there is growing interest in understanding the dynamics of multiscale systems. In this setting the models are typically heuristic in nature, i.e. not derived from first principles, and in many cases the nonlinearities are not presented in an analytic closed form, for example they may be taken as the output of complicated or even black box computer code. However, we assume that the actual system of interest can be usefully modeled by a deterministic dynamical system defined as follows.

DEFINITION 1.1. A *dynamical system* on a topological space $X$ is a continuous map $\varphi : \mathbb{T}^+ \times X \to X$ that satisfies the following two properties:

(i) $\varphi(0, x) = x$ for all $x \in X$, and
(ii) $\varphi(t, \varphi(s, x)) = \varphi(t + s, x)$ for all $s, t \in \mathbb{T}^+$ and all $x \in X$,

where $\mathbb{T}$ denotes the time domain, which is either $\mathbb{Z}$ or $\mathbb{R}$, and $\mathbb{T}^+ := \{t \in \mathbb{T} \mid t \geq 0\}$.

Throughout this paper we make the following assumption

*$X$ is a compact metric space.*

Note that we are *not* assuming that the system is invertible with respect to time. This leads to subtle differences with respect to the more standard theory based on invertible diffeomorphisms or flows (see **[15]**). For example, a set $S \subset X$ is *invariant* under $\varphi$ if

$$\varphi(t, S) = S \quad \forall t \in \mathbb{T}^+.$$

The fact that this condition is restricted to $t \geq 0$, and that $\varphi$ is not assumed to be invertible, makes it important to distinguish between variants of this concept. In Section 2.2 we introduce *forward-backward invariant sets* and *strongly invariant sets*, which are both equivalent to invariant sets for invertible systems. Given a dynamical system $\varphi$ we denote the set of invariant sets by $\mathsf{Invset}(\varphi)$. Moreover, for $U \subset X$ we denote the *maximal invariant set in U* by

$$\mathrm{Inv}(U, \varphi) = \bigcup \big\{ S \subset U \mid S \in \mathsf{Invset}(\varphi) \big\}.$$

For most nonlinear systems that arise in applications, explicit descriptions of the dynamics are usually only obtained via numerical simulations. One of our goals is the development of efficient computational methods which provide mathematically rigorous statements about the structure of the dynamics. As is indicated above, for the applications we have in mind we cannot assume that we have an explicit expression for $\varphi$, but merely an approximation with bounds. Thus our description of the dynamics must be robust with respect to perturbations or errors.

As indicated in the opening sentence of this introduction, it is our belief that understanding the structure of attractors provides a robust computable means of describing the global dynamics. Recall that an *attractor* $A \subset X$ under $\varphi$ can be characterized as a set for which there exists a compact neighborhood $U$ of $A$ such that $A = \omega(U,\varphi) \subset \mathrm{int}\,(U)$, i.e. $A$ is the omega-limit set of $U$ and lies in the interior of $U$. Similarly, a compact set $U \subset X$ is an *attracting neighborhood* under $\varphi$ if $\omega(U,\varphi) \subset \mathrm{int}\,(U)$. The set of attractors and attracting neighborhoods under $\varphi$ are denoted by $\mathsf{Att}(\varphi)$ and $\mathsf{ANbhd}(\varphi)$, respectively. By Proposition 4.3, $\mathsf{Att}(\varphi)$ has the algebraic structure of a bounded distributive lattice. In the case of invertible dynamical systems this has been previously noted in [**14**], however, to accommodate noninvertible systems the operations defining this structure must be modified, as noted in [**1**]. Similarly, Proposition 4.1 guarantees that $\mathsf{ANbhd}(\varphi)$ is a bounded, distributive lattice and by Proposition 3.1, $\mathrm{Inv}(\cdot,\varphi)\colon \mathsf{ANbhd}(\varphi) \to \mathsf{Att}(\varphi)$ is a lattice epimorphism.

A fundamental observation from the computational perspective is that the surjective map $\mathrm{Inv}(\cdot,\varphi)\colon \mathsf{ANbhd}(\varphi) \to \mathsf{Att}(\varphi)$ is stable with respect to perturbations (Proposition 3.22). In contrast, the surjective mapping $\mathrm{Inv}(\cdot,\varphi)\colon 2^X \to \mathsf{Invset}(\varphi)$ applied to all subsets of $X$ typically has no stability properties, i.e. small changes in $U$ can lead to large changes in $\mathrm{Inv}(U,\varphi)$. Moreover, on this level $\mathrm{Inv}$ carries little algebraic structure of invariant sets (it is not a lattice homomorphism), which is in contrast to attractors, as is described below.

Recall that given an attractor $A \in \mathsf{Att}(\varphi)$, its dual *repeller* is defined by $A^* := \{x \in X \mid \omega(x,\varphi) \cap A = \varnothing\}$. Attractor-repeller pairs are the fundamental building blocks for ordering invariant sets in a dynamical system (Theorem 3.19) and lead to Conley's fundamental theorem of dynamical systems [**15**, Theorem IX.1.1]. Repellers can also be characterized in terms of repelling neighborhoods and alpha-limit sets. $\mathsf{Rep}(\varphi)$, the set of repellers under $\varphi$, is a bounded, distributive lattice (Proposition 4.4) as is $\mathsf{RNbhd}(\varphi)$, the set of repelling neighborhoods (Proposition 4.2). While there are lattice anti-homomorphisms between $\mathsf{Att}(\varphi)$ and $\mathsf{Rep}(\varphi)$, and between $\mathsf{ANbhd}(\varphi)$ and $\mathsf{RNbhd}(\varphi)$, the lack of time invertibility of $\varphi$ breaks the complete symmetry between these concepts that is found in the invertible setting. An important example of this is that the lattice epimorphism from repelling neighborhoods to repellers is $\mathrm{Inv}^+(\cdot,\varphi)\colon \mathsf{RNbhd}(\varphi) \to \mathsf{Rep}(\varphi)$ defined by

$$\mathrm{Inv}^+(U,\varphi) := \bigcup \{S \subset U \mid S \in \mathsf{Invset}^+(\varphi)\}$$

where $\mathsf{Invset}^+(\varphi)$ denotes the collection of all forward invariant sets of $\varphi$ as described in Section 2.2.

As is demonstrated in [**10**], given a finite discretization of $X$ and finite approximation of $\varphi$, elements of $\mathsf{ANbhd}(\varphi)$ and $\mathsf{RNbhd}(\varphi)$ can be computed. Part of the motivation for this paper is our belief that useful extensions and proper analysis of these types of computational methods requires a greater understanding of the above mentioned algebraic structures. In particular, it is necessary to clearly distinguish the algebraic operations (which can be handled by the computer) from topological arguments associated with the continuous nature of $X$ and $\varphi$. This separation is not well delineated in the current literature. As a result, the full nature of the algebraic characterization is underutilized; the simplest example being the typical choice a total ordering of a lattice of attractors instead of using the natural partial order.

A key step in this process of delineation is the commutative diagram

$$\begin{array}{ccc} \mathsf{ANbhd}(\varphi) & \overset{c}{\longleftrightarrow} & \mathsf{RNbhd}(\varphi) \\ \omega(\cdot,\varphi)=\mathrm{Inv}(\cdot,\varphi)\Big\downarrow & & \Big\downarrow \alpha(\cdot,\varphi)=\mathrm{Inv}^+(\cdot,\varphi) \\ \mathsf{Att}(\varphi) & \overset{*}{\longleftrightarrow} & \mathsf{Rep}(\varphi) \end{array} \tag{1}$$

where both $^c$ (Proposition 4.6) and $^*$ (Proposition 4.7) are involutions. The proof that this diagram commutes (see Section 4) is the culmination of the topological arguments in Section 3.

The results discussed up to this point have appeared to varying degrees scattered throughout the literature. We hope the following points have been made clear. Attractor-repeller pairs provide a theoretical framework in which to understand the global order structure on nonlinear dynamics. Furthermore, the structure of attractors and repellers can be framed in the algebraic language of lattices, and this algebraic language allows for the transparent development of computer algorithms. In addition, through the maps $\mathrm{Inv}(\cdot,\varphi)$ and $\mathrm{Inv}^+(\cdot,\varphi)$ we can pass from lattices of attracting and repelling neighborhoods, for which there exists efficient computational algorithms, to lattices of attractors and repellers that are not directly computable in general. Finally, these maps are robust with respect to either numerical or experimental perturbation. This strongly suggests that this approach provides a framework in which to develop a new computational theory of nonlinear dynamics, and the first steps in this direction have been taken **[2, 5, 4, 12]**.

For any given dynamical system $\varphi$, $\mathsf{Att}(\varphi)$ has at most countable elements while $\mathsf{ANbhd}(\varphi)$ typically contains uncountably many elements. In general, given $\mathsf{A} \in \mathsf{Att}(\varphi)$, one expects $\mathrm{Inv}^{-1}(\mathsf{A},\varphi)$ to contain uncountably many elements. Of course, for any given approximation scheme and calculation, the collection of attracting neighborhoods that can be explicitly computed is finite, and hence the associated sublattice of attractors is finite. Our goal is to develop a computational theory, which raises the following fundamental question: *Given a finite sublattice of attractors, does there exists a finite sublattice of attracting neighborhoods such that* $\mathrm{Inv}(\cdot,\varphi)$ *produces a lattice isomorphism?* The following theorem answers this question in the affirmative.

THEOREM 1.2. *Let $i$ denote the inclusion map. (i) For every finite sublattice* $\mathsf{A} \subset \mathsf{Att}(\varphi)$*, there exists a lattice monomorphism $k$ such that the following diagram*

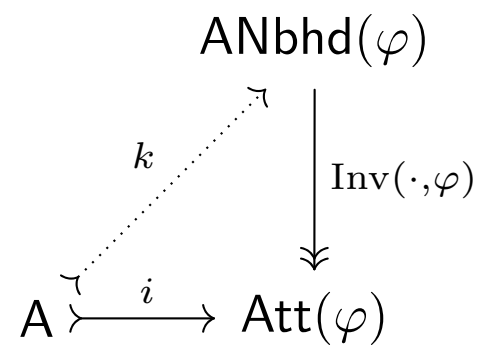

*commutes.*
*(ii) For every finite sublattice* $\mathsf{R} \subset \mathsf{Rep}(X)$*, there exists a lattice monomorphism $k$ such that the following diagram*

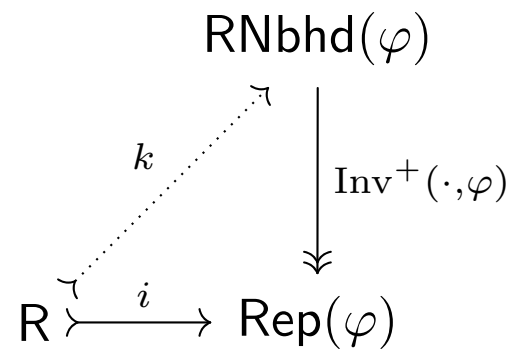

*commutes.*

This theorem implies that there is no fundamental algebraic obstruction to identifying any finite collection of attractors via attracting neighborhoods, which are computable objects [**10**]. Another consequence is an alternative proof of the existence of a index filtration to that presented in [**8**].

It is worth thinking about Theorem 1.2 from a more categorical perspective. Recall that a lattice $\mathsf{H}$ is *projective* if given any diagram of lattice homomorphisms

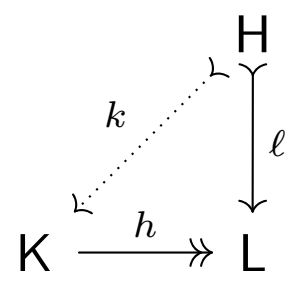

for which the horizontal arrow is an epimorphism, there exists a lattice homomorphism $k$ that makes the diagram commute. Projective, bounded, distributive lattices have been characterized [**3**], and the fact that $\mathsf{Att}(\varphi)$ can have the following lattice

structure

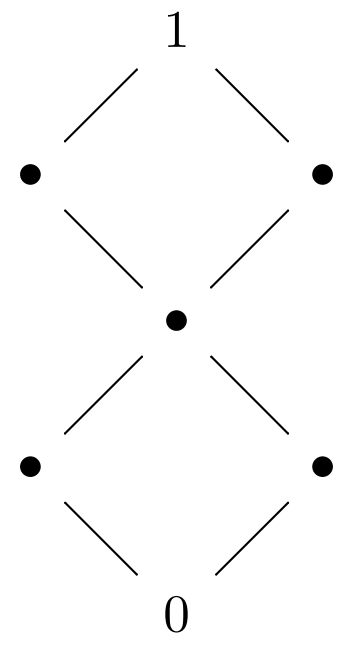

implies that in general a finite sublattice of $\mathsf{Att}(\varphi)$ is not projective.

With this in mind our approach to proving Theorem 1.2 is to introduce the concept of a *spacious* lattice homomorphism (Definition 5.8), and we prove that if in the above diagram $h$ is spacious, then the desired homomorphism $k$ exists (Theorem 5.13). We show that $\mathrm{Inv}^+(\cdot,\varphi)$ is spacious (Proposition 5.12) even though in general $\mathrm{Inv}(\cdot,\varphi)$ is not (Example 5.11). An important benefit of this approach, that is exploited in future work **[11]**, is that if we develop numerical methods under which $\mathrm{Inv}^+(\cdot,\varphi)$ or $\mathrm{Inv}(\cdot,\varphi)$, restricted to computable neighborhoods, is spacious, then we can guarantee that our numerical methods are capable of capturing the desired algebraic structure of the dynamics.

A version of Theorem 1.2 for invertible systems is proved in Robbin and Salamon **[14]**. Without the invertibility assumption, the proof of this theorem is more subtle. In particular, in the invertible case $\mathrm{Inv}(\cdot,\varphi)\colon \mathsf{ANbhd}(\varphi) \to \mathsf{Att}(\varphi)$ is a spacious homomorphism but not in the general case. Therefore, our approach depends on the duality between the lattices of attractors and repellers. In particular, we make use of properties of $\mathrm{Inv}^+(\cdot,\varphi)$ to prove Theorem 1.2(ii) and then use lattice duality between $\mathsf{Att}(\varphi)$ and $\mathsf{Rep}(\varphi)$ to prove (i).

We conclude this introduction with a brief outline. We begin in Section 2.1 with a brief review of essential ideas from lattice theory. We provide a similar review concerning dynamics in Section 2.2.

In Section 3, we present the basic properties of attractors, repellers, and attractor-repeller pairs in the context of potentially noninvertible systems. This culminates in Theorem 3.19 on attractor-repeller pair decompositions that formalizes the idea that attractor-repeller pairs are the fundamental building blocks for ordering invariant sets in a dynamical system. As indicated above one goal is to provide a clear demarcation between the topological arguments and the algebraic computations. Thus this section, which consists of pure point set topological arguments, is followed by Section 4, which contains a demonstration that attractors and repellers have the algebraic structure of bounded, distributive lattices.

A reader whose primary interest is in Theorem 1.2 could accept the results of Sections 3 and 4 and proceed directly to Section 5. In this section we present a proof of this thorem by first proving a general result, Theorem 5.13, concerning the lifting of lattice homomorphisms over bounded, distributive lattices. Then we apply this theorem to the lifting of a finite sublattice of attractors (or repellers) to a lattice of attracting (or repelling) neighborhoods. Theorem 5.13 is of interest in own right and should be useful in the study of the lattice of Lyapunov functions, as in [**14**], and in the context of combinatorial dynamical systems used in computational methods, as mentioned above. These applications will be the subject of future work.

**Acknowledgement.** The first author is partially supported by NSF grant NFS-DMS-0914995, the second author is partially supported by NSF grants NSF-DMS-0835621, 0915019, 1125174, 1248071, and contracts from AFOSR and DARPA. The present work is part of the third authors activities within CAST, a Research Network Program of the European Science Foundation ESF.

## 2. Background

In this section, we summarize the main elements of lattice theory that will be used. Then we review elementary results and introduce notation concerning invariant sets. The results can be viewed as either special cases of the more general framework presented in [**1**] or natural extensions of the standard theory based on invertible dynamics [**15**], and hence we do not provide proofs. In preparation for the work of the later sections we also state and prove elementary lattice properties of these invariant sets. The proofs are fairly elementary, but are included since we are unaware of any single reference for all the results presented.

**2.1. Lattices.** A *lattice* is a set $\mathsf{L}$ with the binary operations $\vee, \wedge : \mathsf{L} \times \mathsf{L} \to \mathsf{L}$ satisfying the following axioms:

(i) (idempotent) $a \wedge a = a \vee a = a$ for all $a \in \mathsf{L}$,
(ii) (commutative) $a \wedge b = b \wedge a$ and $a \vee b = b \vee a$ for all $a, b \in \mathsf{L}$,
(iii) (associative) $a \wedge (b \wedge c) = (a \wedge b) \wedge c$ and $a \vee (b \vee c) = (a \vee b) \vee c$ for all $a, b, c \in \mathsf{L}$,
(iv) (absorption) $a \wedge (a \vee b) = a \vee (a \wedge b) = a$ for all $a, b \in \mathsf{L}$.

A *distributive lattice* satisfies the additional axiom

(v) (distributive) $a \wedge (b \vee c) = (a \wedge b) \vee (a \wedge c)$ and $a \vee (b \wedge c) = (a \vee b) \wedge (a \vee c)$ for all $a, b, c \in \mathsf{L}$.

For distributivity, if one of the two conditions in (v) is satisfied, then so is the other. A lattice is *bounded* if there exist *neutral* elements 0 and 1 with property that

(vi) $0 \wedge a = 0, \quad 0 \vee a = a, \quad 1 \wedge a = a,$ and $1 \vee a = 1$ for all $a \in \mathsf{L}$.

A subset $\mathsf{K} \subset \mathsf{L}$ is called a *sublattice* of $\mathsf{L}$, if $a, b \in \mathsf{K}$ implies that $a \vee b \in \mathsf{K}$ and $a \wedge b \in \mathsf{K}$. For sublattices of bounded lattices we impose the additional condition that $0, 1 \in \mathsf{K}$.

Let $\mathsf{K}$ and $\mathsf{L}$ be lattices. A function $h : \mathsf{K} \to \mathsf{L}$ is a *lattice homomorphism* if

$$h(a \wedge b) = h(a) \wedge h(b) \text{ and } h(a \vee b) = h(a) \vee h(b),$$

and a *lattice anti-homomorphism* if $\vee$ is replaced by $\wedge$ and vice versa. A lattice homomorphism between bounded lattices is a lattice homomorphism with the additional property that $h(0) = 0$ and $h(1) = 1$. Bounded, distributive lattices together with the above described morphisms form a category denoted by $\mathsf{BDLat}$. Summarizing: *In the category of bounded, distributive lattices, sublattices contain* $0$ *and* $1$*, and the lattice homomorphisms preserve* $0$ *and* $1$.

In some of the structures we encounter in this paper inversion relations are satisfied. In the context of lattices this is often described by Boolean algebras. A *Boolean algebra* is a bounded, distributive lattice $\mathsf{B}$ with the additional complementation relation $a \mapsto a^c \in \mathsf{B}$ which satisfies the axiom:

(vii) $a \wedge a^c = 0$ and $a \vee a^c = 1$ for all $a \in \mathsf{B}$.

A Boolean homomorphism between Boolean algebras is a lattice homomorphism which preserves complements. We remark that if $h : \mathsf{B} \to \mathsf{C}$ is a lattice homomorphism between Boolean algebras, then $f$ is automatically a Boolean homomorphism (see [**7**, Lemma 4.17]). Boolean algebras with Boolean homomorphisms form a category, denoted by $\mathsf{Bool}$.

In this context, there is the following classical result: *Every bounded, distributive lattice can be embedded in a Boolean algebra*, see [**9**, Theorem 153]. However, we need the following stronger result that establishes this procedure as a functor.

THEOREM 2.1. *Given a bounded, distributive lattice* $\mathsf{L}$*, then there is a unique (up to isomorphism) Boolean algebra* $\mathsf{B}(\mathsf{L})$ *and a lattice homomorphism* $j : L \to \mathsf{B}(\mathsf{L})$ *with the property that for every homomorphism* $g$ *from* $\mathsf{L}$ *to a Boolean algebra* $\mathsf{C}$ *there exists a unique lattice homomorphism* $\mathsf{B}(g) : \mathsf{B}(\mathsf{L}) \to \mathsf{C}$ *such that* $\mathsf{B}(g) \circ j = g$.

The Boolean algebra $\mathsf{B}(\mathsf{L})$ in the above theorem is the *Booleanization of* $\mathsf{L}$, and the theorem implies that Booleanization is a covariant functor from the category of bounded, distributive lattices $\mathsf{BDLat}$ to the category of Boolean algebras $\mathsf{Bool}$, see [**17**, Definition 9.5.5] and [**13**, Corollary 20.11].

Our approach makes significant use of Birkhoff's Representation Theorem, which provides a deep relation between finite distributive lattices and posets. To state this theorem requires the introduction of several concepts that allow us to move back and forth between posets and lattices.

A *poset* $(\mathsf{P}, \leq)$ is a set $\mathsf{P}$ with a binary relation $\leq$, called a *partial order*, which satisfies the following axioms:

(i) (reflexivity) $p \leq p$ for all $p \in \mathsf{P}$,
(ii) (anti-symmetry) $p \leq q$ and $q \leq p$ implies $p = q$,
(iii) (transitivity) if $p \leq q$ and $q \leq r$, then $p \leq r$.

Let $\mathsf{P}$ and $\mathsf{Q}$ be posets. A mapping $f : \mathsf{P} \to \mathsf{Q}$ is called *order-preserving* if $f(p) \leq f(q)$ for all $p \leq q$. A mapping is an *order-embedding* if $f(p) \leq f(q)$ if and only if $p \leq q$. Posets together with order-preserving mappings form a category denoted by $\mathsf{Poset}$.

A lattice $\mathsf{L}$ has a naturally induced partial order as follows. Given $a, b \in \mathsf{L}$ define

$$a \leq b \quad \Leftrightarrow \quad a \wedge b = a \quad \Leftrightarrow \quad a \wedge b^c = 0, \tag{2}$$

where the latter relation only makes sense in the setting of a Boolean algebra.

Given an element $p \in \mathsf{P}$, the *down-set* and *up-set* of $p$ are the sets $\downarrow p = \{q \in \mathsf{P} \mid q \leq p\}$ and $\uparrow p = \{q \in \mathsf{P} \mid p \leq q\}$, respectively. The collection of down-sets of a finite poset $\mathsf{P}$ generates a finite distributive lattice denoted by $\mathsf{O}(\mathsf{P})$ with respect to $\vee = \cup$ and $\wedge = \cap$.

Given any poset $\mathsf{P}$ its *dual*, denoted by $\mathsf{P}^\partial$, is defined to be the poset with order $p \leq q$ in $\mathsf{P}^\partial$ if and only if $q \leq p$ in $\mathsf{P}$. It is left to the reader to check that the *complement* function

$$\begin{aligned} {}^c \colon \mathsf{O}(\mathsf{P}) \quad &\to \quad \mathsf{O}(\mathsf{P}^\partial) \\ \alpha \quad &\mapsto \quad \alpha^c := \mathsf{P} \setminus \alpha \end{aligned}$$

is a involutive lattice anti-morphism.

Given a lattice $\mathsf{L}$, an element $c \in \mathsf{L}$ is *join-irreducible* if

(a) $c \neq 0$ and
(b) $c = a \vee b$ implies $c = a$ or $c = b$ for all $a, b \in \mathsf{L}$.

The set of join-irreducible elements in $\mathsf{L}$ is denoted by $\mathsf{J}(\mathsf{L})$. Observe that $c$ is join-irreducible if and only if there exists a unique element $a \in \mathsf{L}$ such that $a < c$, and hence we can define the predecessor map $\overleftarrow{\ } \colon \mathsf{J}(\mathsf{L}) \to \mathsf{L}$. The set $\mathsf{J}(\mathsf{L})$ is a poset as a subset of $\mathsf{L}$.

Let $\mathsf{FDLat}$ denote the small category of finite distributive lattices, whose homomorphisms map $0$ and $1$ to itself according to our conventions, and let $\mathsf{FPoset}$ denote the small category of finite posets, whose morphisms are order-preserving mappings. For a finite distributive lattice $\mathsf{L}$ and finite poset $\mathsf{P}$ the maps

$$\mathsf{L} \xrightarrow{\mathsf{J}} \mathsf{J}(\mathsf{L}) \quad \text{and} \quad \mathsf{P} \xrightarrow{\mathsf{O}} \mathsf{O}(\mathsf{P})$$

define contravariant functors from $\mathsf{FDLat}$ to $\mathsf{FPoset}$ and from $\mathsf{FPoset}$ to $\mathsf{FDLat}$ respectively.

Finally, define the mapping

$$\begin{aligned} \downarrow^\vee \colon \mathsf{L} \quad &\to \quad \mathsf{O}(\mathsf{J}(\mathsf{L})) \\ a \quad &\mapsto \quad \downarrow^\vee a = \{b \in \mathsf{J}(\mathsf{L}) \mid b \leq a\} = \downarrow a \cap \mathsf{J}(\mathsf{L}). \end{aligned}$$

THEOREM 2.2 (Birkhoff's Representation Theorem). *Let* $\mathsf{L}$ *be a finite distributive lattice and let* $\mathsf{P}$ *be a finite partially ordered set. Then* $\downarrow^{\vee}\colon \mathsf{L} \to \mathsf{O}\big(\mathsf{J}(\mathsf{L})\big)$ *is a lattice isomorphism, and* $\downarrow\colon \mathsf{P} \to \mathsf{J}\big(\mathsf{O}(\mathsf{P})\big)$ *is an order isomorphism.*

For a proof see **[7, 16]**. We can use this representation to recast Theorem 2.1.

PROPOSITION 2.3. *Let* $\mathsf{P}$ *be a finite poset, and let* $f\colon \mathsf{O}(\mathsf{P}) \to \mathsf{L}$ *be a lattice homomorphism, where* $\mathsf{L}$ *is a bounded, distributive lattice with Booleanization* $\mathsf{B}(\mathsf{L})$ *and lattice-embedding* $j\colon \mathsf{L} \to \mathsf{B}(\mathsf{L})$*. Then there exists a unique homomorphism* $\mathsf{B}(f)\colon 2^{\mathsf{P}} = \mathsf{B}(\mathsf{O}(\mathsf{P})) \to \mathsf{B}(\mathsf{L})$ *such that for every* $\alpha \in \mathsf{O}(\mathsf{P}) \subset 2^{\mathsf{P}}$ *we have* $j \circ f(\alpha) = \mathsf{B}(f)(\alpha)$*. In particular, for every* $\alpha \in 2^{\mathsf{P}}$ *we have*

$$\mathsf{B}(f)(\alpha) = \bigvee_{p\in\alpha} \mathsf{B}(f)(\{p\}). \tag{3}$$

*Moreover,*

$$\mathsf{B}(f)(\{p\}) \wedge \mathsf{B}(f)(\{p'\}) = 0 \tag{4}$$

*for all* $p, p' \in \mathsf{P}$ *with* $p \neq p'$*.*

PROOF. The result follows immediately from Theorem 2.1 by observing that since $\mathsf{P}$ is finite, we can take $\mathsf{B}(\mathsf{O}(\mathsf{P})) = 2^{\mathsf{P}}$ where the corresponding embedding $i\colon \mathsf{O}(\mathsf{P}) \to 2^{\mathsf{P}}$ is the inclusion map. ∎

**2.2. Invariant sets.** Since we assume $\varphi$ to be defined and single-valued only for positive times, $\varphi$ does not necessarily extend to a function satisfying Definition 1.1(ii) for negative times. However, for $t < 0$, the sets

$$\varphi(t, x) := \{y \in X \mid \varphi(-t, y) = x\}$$

extend $\varphi$ to a (possibly) multivalued map for negative times, but note that $\varphi(t, x)$ may be the empty set in this case. The semigroup property implies that $\varphi(t, X) = X$ for some $t > 0$ if and only if $\varphi(t, X) = X$ for all $t \in \mathbb{T}^+$. In this case $\varphi$ defines a *surjective dynamical system*, for which $\varphi(t, S) \neq \varnothing$ for all $\varnothing \neq S \subset X$ and $t \in \mathbb{T}$. Note that a dynamical system is surjective if and only if $X$ is an invariant set.

For $x \in X$, consider a continuous function $\gamma_x\colon \mathbb{T} \to X$, satisfying $\gamma_x(0) = x$, and $\gamma_x(t+s) = \varphi(s, \gamma_x(t))$ for all $t \in \mathbb{T}$ and all $s \in \mathbb{T}^+$. Its image, also denoted by $\gamma_x$, is called a *complete orbit* through $x$. A set $S$ is invariant if and only if for each $x \in S$ there exists a complete orbit $\gamma_x \subset S$. The restriction to $t \leq 0$ gives the *backward orbit* $\gamma_x^-$ and the restriction to $t \geq 0$ gives the *forward orbit* $\gamma_x^+$. The backward image of a set $U$ is $\Gamma^-(U) = \bigcup_{t\leq 0} \varphi(t, U)$ and the forward image is $\Gamma^+(U) = \bigcup_{t\geq 0} \varphi(t, U)$.

Recall that a set $S \subset X$ is invariant if

$$\varphi(t, S) = S \quad \forall t \in \mathbb{T}^+.$$

A set $S \subset X$ is *forward invariant* if

$$\varphi(t, S) \subset S \quad \forall t \in \mathbb{T}^+.$$

The collection of forward invariant sets is denoted by $\mathsf{Invset}^+(\varphi)$. Observe that $\mathsf{Invset}(\varphi) \subset \mathsf{Invset}^+(\varphi)$. Similarly, a set $S$ is *backward invariant* if

$$\varphi(t, S) \subset S \quad \forall t \in \mathbb{T}^-.$$

and backward invariant sets are denoted by $\mathsf{Invset}^-(\varphi)$. Sets that are both forward and backward invariant are called *forward-backward invariant* sets and denoted by $\mathsf{Invset}^\pm(\varphi) = \mathsf{Invset}^+(\varphi) \cap \mathsf{Invset}^-(\varphi)$. Observe that forward invariance implies that $S \subset \varphi(-t, \varphi(t, S)) \subset \varphi(-t, S)$, for $t \geq 0$ and backward invariance implies that $\varphi(t, S) \subset S$, for $t \leq 0$, thus elements of $\mathsf{Invset}^\pm(\varphi)$ can be characterized by:

$$\varphi(t, S) = S \quad \forall t \in \mathbb{T}^-$$

In particular, the phase space $X \in \mathsf{Invset}^\pm(\varphi)$. Finally a set $S \subset X$ is called *strongly invariant* if

$$\varphi(t, S) = S \quad \forall t \in \mathbb{T}.$$

The strongly invariant sets are denoted by $\mathsf{SInvset}(\varphi) = \mathsf{Invset}(\varphi) \cap \mathsf{Invset}^\pm(\varphi)$. If $\varphi$ is surjective, then $\mathsf{SInvset}(\varphi) = \mathsf{Invset}^\pm(\varphi)$. For invertible systems $\mathsf{Invset}(\varphi) = \mathsf{Invset}^\pm(\varphi) = \mathsf{SInvset}(\varphi)$.

EXAMPLE 2.4. Consider $\varphi \colon \mathbb{Z}^+ \times [-1, 1] \to [-1, 1]$ generated by

$$f(x) = \begin{cases} 0 & \text{for } x \in [-1, 0] \\ \frac{5}{2}x(1-x) & \text{for } x \in [0, 1]. \end{cases}$$

Observe that $K = \{\frac{3}{5}\} \in \mathsf{Invset}(\varphi)$, but $f^{-1}(K) = \{\frac{2}{5}, \frac{3}{5}\}$ and hence $K \notin \mathsf{SInvset}(\varphi)$. Thus $\mathsf{Invset}(\varphi) \not\subset \mathsf{SInvset}(\varphi) \subset \mathsf{Invset}^-(\varphi)$.

We leave it to the reader to check that a dynamical system $\varphi$ can be restricted to $S \in \mathsf{Invset}^+(\varphi)$, that is,

$$\begin{array}{rcl} \varphi|_S \colon \mathbb{T}^+ \times S & \to & S \\ (t, x) & \mapsto & \varphi(t, x) \end{array}$$

is a dynamical system. The backward extension is given by $\varphi|_S(t, x) = \varphi(t, x) \cap S$ for all $t < 0$. If $S \in \mathsf{Invset}(\varphi)$, then $\varphi|_S$ is a surjective dynamical system on $S$.

The sets $\mathsf{Invset}^+(\varphi)$ and $\mathsf{Invset}^-(\varphi)$ are bounded, distributive lattices with respect to the binary operations $\vee = \cup$ and $\wedge = \cap$. The neutral elements are $0 = \varnothing$ and $1 = X$.

PROPOSITION 2.5. *The function*

$$\begin{array}{rcl} {}^c \colon \mathsf{Invset}^+(\varphi) & \to & \mathsf{Invset}^-(\varphi) \\ S & \mapsto & S^c \end{array}$$

*is an involutive lattice anti-isomorphism.*

PROOF. De Morgan's laws imply that on the level of sets the complement map is a lattice anti-homomorphism.

If $S \in \mathsf{Invset}^+(\varphi)$, then $\varphi(t,S) \subset S$, for all $t \geq 0$. Therefore $S \subset \varphi(-t,\varphi(t,S)) \subset \varphi(-t,S)$, and hence $S \subset \varphi(t,S)$ for all $t \leq 0$. Since $\varphi(t,X) = X$ for all $t \leq 0$, we have $S^c \supset \varphi(t,S)^c = \varphi(t,S^c)$ for all $t \leq 0$, which proves that $S^c \in \mathsf{Invset}^-(\varphi)$.

If $S \in \mathsf{Invset}^-(\varphi)$, then $\varphi(t,S) \subset S$ for all $t \leq 0$. Moreover, for $t \leq 0$, $X = \varphi(t,X) = \varphi(t,S\cup S^c) = \varphi(t,S) \cup \varphi(t,S^c)$, from which we derive that $S^c \subset \varphi(t,S^c)$ for all $t \leq 0$. Furthermore, $\varphi(-t,S^c) \subset \varphi(-t,\varphi(t,S^c)) \subset S^c$ for all $t \leq 0$, and thus $\varphi(t,S^c) \subset S^c$ for $t \geq 0$, which proves that $S^c \in \mathsf{Invset}^+(\varphi)$. ∎

COROLLARY 2.6. $\mathsf{Invset}^\pm(\varphi)$ *is a Boolean algebra.*

Next we establish that $\mathsf{Invset}(\varphi)$ is also a lattice. Through out this paper, the dynamical system $\varphi$ is considered fixed and thus for the sake of notational simplicity we let

$$\mathrm{Inv}(\cdot) = \mathrm{Inv}(\cdot,\varphi) \quad \text{and} \quad \mathrm{Inv}^+(\cdot) = \mathrm{Inv}^+(\cdot,\varphi).$$

Observe that $\mathsf{Invset}(\varphi)$ (as well as $\mathsf{Invset}^+(\varphi), \mathsf{Invset}^-(\varphi)$, and $\mathsf{Invset}^\pm(\varphi)$) is a poset with respect to inclusion and $\mathrm{Inv}\colon 2^X \to \mathsf{Invset}(\varphi)$ is order-preserving. We leave it to the reader to check that if $S,S' \in \mathsf{Invset}(\varphi)$, then $S\cup S' \in \mathsf{Invset}(\varphi)$ and $S \cup S' = \sup(S,S')$. In general, $S \cap S' \in \mathsf{Invset}^+(\varphi)$ and hence $\mathsf{Invset}(\varphi)$ is *not* necessarily closed under intersection. However, as indicated below, every pair $S,S' \in \mathsf{Invset}(\varphi)$ has a greatest lower bound.

LEMMA 2.7. *For every pair* $S,S' \in \mathsf{Invset}^-(\varphi)$ *or* $S,S' \in \mathsf{Invset}^+(\varphi)$*,*

$$\begin{aligned} \mathrm{Inv}(S\cup S') &= \mathrm{Inv}(S) \cup \mathrm{Inv}(S'); \\ \mathrm{Inv}(S\cap S') &= \mathrm{Inv}\big(\mathrm{Inv}(S)\cap \mathrm{Inv}(S')\big). \end{aligned}$$

PROOF. We prove the lemma for forward invariant sets. The arguments are similar for backward invariant sets. Let $S,S' \in \mathsf{Invset}^+(\varphi)$, then $\mathrm{Inv}(S)\cup\mathrm{Inv}(S') \subset \mathrm{Inv}(S\cup S')$. To establish the reverse inclusion we argue as follows. Let $x \in \mathrm{Inv}(S\cup S')$ and let $\gamma_x \subset \mathrm{Inv}(S\cup S')$ be a complete orbit. Since $S,S' \in \mathsf{Invset}^+(\varphi)$, it follows that for all $y \in \gamma_x$ we have that $\gamma_y^+ \subset S$ when $y \in S$, and $\gamma_y^+ \subset S'$ when $y \in S'$. This implies that $y \in S$ or $y \in S'$ for all $y \in \gamma_x$, and thus $\gamma_x \in \mathrm{Inv}(S)$ or $\gamma_x \in \mathrm{Inv}(S')$. Consequently, $\mathrm{Inv}(S\cup S') \subset \mathrm{Inv}(S)\cup\mathrm{Inv}(S')$.

As for intersections we have that $\mathrm{Inv}(S)\cap\mathrm{Inv}(S') \subset S\cap S'$, and therefore $\mathrm{Inv}\big(\mathrm{Inv}(S)\cap\mathrm{Inv}(S')\big) \subset \mathrm{Inv}(S\cap S')$. On the other hand $\mathrm{Inv}(S\cap S') \subset \mathrm{Inv}(S)$, and $\mathrm{Inv}(S\cap S') \subset \mathrm{Inv}(S')$, which implies $\mathrm{Inv}(S\cap S') \subset \mathrm{Inv}(S)\cap\mathrm{Inv}(S')$ and $\mathrm{Inv}(S\cap S') \subset \mathrm{Inv}\big(\mathrm{Inv}(S)\cap\mathrm{Inv}(S')\big)$. ∎

PROPOSITION 2.8. *With the binary operations* $\vee = \sup = \cup$ *and* $\wedge$ *defined by*

$$S \wedge S' = \inf(S,S') = \mathrm{Inv}(S\cap S'),$$

$\mathsf{Invset}(\varphi)$ *is a bounded distributive lattice. The neutral elements are* $0 = \varnothing$ *and* $1 = \mathrm{Inv}(X)$.

PROOF. We prove the distributive property. For $S, S', S'' \in \mathsf{Invset}(X)$, then sets $S \cap S'$ and $S \cap S''$ are forward invariant. Then by Lemma 2.7

$$\begin{aligned}(S \wedge S') \vee (S \wedge S'') &= \mathrm{Inv}(S \cap S') \cup \mathrm{Inv}(S \cap S'') \\ &= \mathrm{Inv}((S \cap S') \cup (S \cap S'')) \\ &= \mathrm{Inv}(S \cap (S' \cup S'')) = \mathrm{Inv}(S) \wedge \mathrm{Inv}(S' \cup S'') \\ &= S \wedge (S' \vee S''),\end{aligned}$$

which completes the proof. ∎

Lemma 2.7 implies that $\mathrm{Inv}\colon \mathsf{Invset}^-(\varphi) \to \mathsf{Invset}(\varphi)$ and $\mathrm{Inv}\colon \mathsf{Invset}^+(\varphi) \to \mathsf{Invset}(\varphi)$ are lattice homomorphisms.

$\mathsf{SInvset}(\varphi)$ is a lattice with respect to intersection and union, and a sublattice of $\mathsf{Invset}(\varphi)$. The following lemma establishes an important property relating to the algebra of invariant sets (see Proposition 3.13).

LEMMA 2.9. *Let* $S \in \mathsf{Invset}^{\pm}(\varphi)$ *and* $S' \in \mathsf{Invset}(\varphi)$, *then* $S \cap S' \in \mathsf{Invset}(\varphi)$.

PROOF. Since $S$ and $S'$ are both forward invariant, it follows that $\varphi(t, S \cap S') \subset S \cap S'$ for all $t \geq 0$. Since $S \in \mathsf{Invset}^{\pm}(\varphi)$, it holds that $\varphi(-t, S) = S$ for all $t \geq 0$. The invariance of $S'$ implies that $\varphi|_{S'}$ is surjective, and therefore for each $x \in S'$ we have $\varphi(-t, x) \cap S' \neq \varnothing$. Combining these facts yields that for each $x \in S \cap S'$ and each $t \geq 0$ there exists a point $y \in S \cap S'$ such that $\varphi(t, y) = x$. This shows that $S \cap S' \subset \varphi(t, S \cap S')$ for all $t \geq 0$. Therefore $\varphi(t, S \cap S') = S \cap S'$ for all $t \geq 0$, which proves the invariance of $S \cap S'$. ∎

LEMMA 2.10. $\mathrm{Inv}^+\colon \mathsf{Invset}^-(\varphi) \to \mathsf{Invset}^{\pm}(\varphi)$.

PROOF. By definition $\mathrm{Inv}^+(U)$ is forward invariant. Since $U$ is backward invariant, then for every $x \in \mathrm{Inv}^+(U)$ we have $\Gamma^-(x) \subset U$. Let $y \in \Gamma^-(x)$, then $\varphi(s, y) = x$ for some $s \geq 0$. Consequently, $\varphi(t + s, y) = \varphi(t, \varphi(s, y)) = \varphi(t, x) \subset U$, and $\varphi(\sigma, y) \in U$ for $0 \leq \sigma \leq s$, and thus $\gamma_y^+ \subset U$. This implies $y \in \mathrm{Inv}^+(U)$ and therefore $\Gamma^-(x) \subset \mathrm{Inv}^+(U)$, which proves that $\varphi(t, \mathrm{Inv}^+(U)) \subset \mathrm{Inv}^+(U)$ for all $t \leq 0$. ∎

The reader can check that $\mathrm{Inv}^+\colon \mathsf{Invset}^-(\varphi) \to \mathsf{Invset}^{\pm}(\varphi)$ is a lattice homomorphism.

For a set $U \subset X$ define

$$\alpha(U) = \bigcap_{t \leq 0} \mathrm{cl}\,\Big(\varphi\big((-\infty, t], U\big)\Big) \quad \text{and} \quad \omega(U) = \bigcap_{t \geq 0} \mathrm{cl}\,\Big(\varphi\big([t, \infty), U\big)\Big),$$

which are called the *alpha-limit and omega-limit sets* of $U$ respectively. For noninvertible dynamical systems, as studied in this paper, there is a lack of symmetry between alpha-limit and omega-limit sets. In the next two propositions we list the important properties of limit sets that are used throughout the remainder of this work.

PROPOSITION 2.11. *Let $U \subset X$ , then $\omega(U)$ satisfies the following list of properties:*

(i) $\omega(U) \in \mathsf{Invset}(\varphi)$ *is compact;*
(ii) $U \neq \varnothing$ *implies* $\omega(U) \neq \varnothing$*;*
(iii) *if* $\varphi(t, U) \subset U$ *for all* $t \geq \tau \geq 0$*, then* $\omega(U) = \mathrm{Inv}\big(\mathrm{cl}\,(U)\big) \subset \mathrm{cl}\,(U)$*;*
(iv) $V \subset U$ *implies* $\omega(V) \subset \omega(U)$*;*
(v) $\omega(U \cup V) = \omega(U) \cup \omega(V)$ *and* $\omega(U \cap V) \subset \omega(U) \cap \omega(V)$*;*
(vi) $\omega(U) = \omega\big(\mathrm{cl}\,(U)\big)$*;*
(vii) *if there exists a backward orbit* $\gamma_x^- \subset U$*, then* $x \in \omega(U)$*;*
(viii) *if* $U \in \mathsf{Invset}(\varphi)$*, then* $\mathrm{cl}\,(U) = \omega(U)$*.*

REMARK 2.12. Proposition 2.11(viii) implies that the closure of an invariant set is invariant and the closure of a forward invariant is forward invariant, i.e. $\mathrm{cl} : \mathsf{Invset} \to \mathsf{Invset}$ and $\mathrm{cl} : \mathsf{Invset}^+ \to \mathsf{Invset}^+$. However, this property does not hold for backward invariant sets, cf. Example 2.14.

PROPOSITION 2.13. *Let $U \subset X$ (with $X$ compact), then $\alpha(U)$ satisfies the following list of properties:*

(i) $\alpha(U) \in \mathsf{Invset}^+(\varphi)$ *is compact;*
(ii) $U \neq \varnothing$ *and* $\varphi$ *surjective implies* $\alpha(U) \neq \varnothing$*;*
(iii) *if* $U \in \mathsf{Invset}^-(\varphi)$*, then* $\alpha(U) \subset \mathrm{cl}\,(U)$*, and if* $\mathrm{cl}\,(U) \in \mathsf{Invset}^-(\varphi)$*, then* $\alpha\big(\mathrm{cl}\,(U)\big) = \mathrm{Inv}^+\big(\mathrm{cl}\,(U)\big)$*;*
(iv) *if* $V \subset U$*, then* $\alpha(V) \subset \alpha(U)$*;*
(v) $\alpha(U \cup V) = \alpha(U) \cup \alpha(V)$ *and* $\alpha(U \cap V) \subset \alpha(U) \cap \alpha(V)$*;*
(vi) *if* $\gamma_x^+ \subset U$*, then* $x \in \alpha(U)$*. In particular* $\mathrm{Inv}^+(U) \subset \alpha(U)$*, and* $\mathrm{Inv}^+(U) = \alpha(U)$ *whenever* $\alpha(U) \subset U$*;*
(vii) *if* $U \in \mathsf{Invset}^-(\varphi)$ *and compact, then* $\alpha(U) \in \mathsf{Invset}^\pm(\varphi)$*. If in addition* $\varphi$ *is surjective, then* $\alpha(U) \in \mathsf{SInvset}(\varphi)$ *and* $\alpha(U) = \mathrm{Inv}(U)$*;*
(viii) $U \in \mathsf{Invset}^+(\varphi)$ *implies* $\mathrm{cl}\,(U) \subset \alpha(U)$*, and if* $U \in \mathsf{Invset}^\pm(\varphi)$*, then* $\mathrm{cl}\,(U) = \alpha(U)$*.*

For a point $x \in X$ and for a backward orbit $\gamma_x^-$ we define the *orbital alpha-limit set* by

$$\alpha_{\mathrm{o}}(\gamma_x^-) = \bigcap_{t \leq 0} \mathrm{cl}\,\big(\gamma_x\big((-\infty, t]\big)\big).$$

EXAMPLE 2.14. Alpha-limit sets need not be invariant. Consider the dynamical system in Example 2.4. Let $U = [-1, 0]$. Then $\alpha(U) = U$, but $f(U) = \{0\}$ and hence $\alpha(U) \notin \mathsf{Invset}(\varphi)$.

Alpha-limit sets need not be nonempty. Consider $W = [-1, 0)$, then $\alpha(W) = \varnothing$. This example also shows that Proposition 2.13(viii) is sharp as $W \in \mathsf{Invset}^{-}(\varphi)$ and $\mathrm{cl}\,(W) \not\subset \alpha(W)$.

Alpha-limit sets need not be backward invariant. Let $z = \frac{3}{5}$. There are two distinct backward orbits: the constant orbit $\gamma_c^-(n) = \frac{3}{5}$ for all $n \in \mathbb{Z}^-$; and the orbit $\gamma^-(n) = f_-^{-n}(\frac{3}{5})$ for $n \in \mathbb{Z}^-$ where $f_-^{-1}(x) := (1 - \sqrt{1-(8x)/5})/2$. Since $\gamma^-(n) \to 0$ as $n \to -\infty$, $\alpha(\frac{3}{5}) = \gamma^- \cup \{0\}$. Then $f^{-1}(\alpha(\frac{3}{5})) = \gamma^- \cup [-1, 0]$, which shows that $\alpha(\frac{3}{5}) \notin \mathsf{Invset}^{-}(\varphi)$.

Now, consider the set $V = \{\frac{3}{5}\}$, which is invariant but not backward invariant. In particular there are two distinct backward orbits from the point $\frac{3}{5}$, namely the constant orbit $\gamma_c^-(n) = \frac{3}{5}$ for all $n \in \mathbb{Z}^-$ and the orbit $\gamma^-(n) = f_-^{-n}(\frac{3}{5})$ for $n \in \mathbb{Z}^-$ where $f_-^{-1}(x) := (1 - \sqrt{1-(8x)/5})/2$. Since $\gamma^-(n) \to 0$ as $n \to -\infty$, $\alpha(V) = \gamma^- \cup \{0\}$. Then $f^{-1}(\alpha(V)) = \gamma^- \cup [-1, 0]$, which shows that $\alpha(V)$ is not backward invariant. Moreover, $\alpha_{\mathrm{o}}(\gamma_x^-) = \{0\} \subsetneq \alpha(\{x\})$ with $x = 3/5$.

To see that compactness is necessary in Proposition 2.13(vii), consider the full, nonconstant orbit $\gamma$ through $\frac{3}{5}$. Observe that $\gamma \in \mathsf{Invset}^{-}(\varphi)$, but $\mathrm{cl}\,(\gamma) = \gamma \cup \{0\} \notin \mathsf{Invset}^{-}(\varphi)$.

PROPOSITION 2.15. *If $x \in U \subset X$, then $\alpha_o(\gamma_x^-)$ is nonempty, compact, and invariant, and $\alpha_{\mathrm{o}}(\gamma_x^-) \subset \alpha(U)$.*

Orbital alpha-limit sets and omega-limit sets can be used to introduce a notion of dual set of a subset $S \subset X$. The *positive dual set* of $S$ is

$$S^{\oplus} := \{x \in X \mid \omega(x) \cap S = \varnothing\}$$

and the *negative dual set* of $S$ is

$$S^{\ominus} := \{x \in X \mid \exists \gamma_x \text{ such that } \alpha_{\mathrm{o}}(\gamma_x^-) \cap S = \varnothing\}.$$

PROPOSITION 2.16. *Let $S \subset X$, then $S^{\oplus}$ is forward-backward invariant, and $S^{\ominus}$ is invariant. Moreover, if $S$ is compact, then invariance of $S$ implies that $S \cap S^{\oplus} = \varnothing$, and forward-backward invariance of $S$ implies that $S \cap S^{\ominus} = \varnothing$.*

## 3. Attractors and Repellers

A *trapping region* is a forward invariant set $U \subset X$ such that

$$\varphi(\tau, \mathrm{cl}\,(U)) \subset \mathrm{int}\,(U) \quad \text{for some } \tau > 0.$$

A set $A \subset X$ is an *attractor* if there exists a trapping region $U$ such that $A = \mathrm{Inv}(U)$.

**3.1. Properties of attractors.** As indicated in the introduction one of ours goals is the commutative diagram (1). To obtain this requires a detailed understanding of the relationship between attractors and their neighborhoods. This is a well studied topic and thus we only state the specific results that we make use of. Proofs of these results can either be found in [**15**, Section IX.9.1] or can be obtained using the similar arguments. We begin by providing an equivalent characterization of attractors.

PROPOSITION 3.1. *If $U$ is a trapping region, then the associated attractor* $\mathrm{Inv}(U) = \mathrm{Inv}(\mathrm{cl}\,(U)) = \omega(U)$ *is compact.*

The following lemmata lead to a proof of Proposition 3.5 (alternatively see [**15**, Proposition IX.9.1.10]), but are included because they are explicitly cited in the proofs provided in late subsections of this section.

LEMMA 3.2. *A neighborhood $U \subset X$ is an attracting neighborhood if and only if there exists a $\tau > 0$ such that $\varphi(t, \mathrm{cl}\,(U)) \subset \mathrm{int}\,(U)$, for all $t \geq \tau$. In addition, if $U$ is an attracting neighborhood, then so are* $\mathrm{cl}\,(U)$ *and* $\mathrm{int}\,(U)$.

LEMMA 3.3. *A neighborhood $U \subset X$ is a trapping region if and only if $U$ is a forward invariant attracting neighborhood.*

LEMMA 3.4. *Let $U \subset X$ be attracting neighborhood for an attractor $A = \omega(U)$, and let $U' \subset X$ be a neighborhood of $A$ such that $A \subset U' \subset \mathrm{cl}\,(U)$. Then also $U'$ is an attracting neighborhood with $A = \omega(U')$.*

The following proposition asserts that the omega-limit set of an attracting neighborhood is an attractor. Therefore we say that *$U$ is an attracting neighborhood for an attractor $A$* if $A = \omega(U) \subset \mathrm{int}\,(U)$.

PROPOSITION 3.5. *A set $A \subset X$ is an attractor if and only if there exists a neighborhood $U$ of $A$ such that $A = \omega(U)$. In particular, $U$ is an attracting neighborhood. Moreover, for every attracting neighborhood $U$ there exists a trapping region $U' \subset U$.*

COROLLARY 3.6. *Let $U \subset X$ be an attracting neighborhood, then* $A = \mathrm{Inv}(U) = \mathrm{Inv}(\mathrm{cl}\,(U)) = \omega(U)$ *is an attractor.*

**3.2. An attractor inside an attractor.** An important property of attractors concerns relative attractors inside an attractor.

PROPOSITION 3.7. *If $A \subset X$ is an attractor for $\varphi$ and $A' \subset A$ is an attractor for $\varphi|_A$, then $A'$ is an attractor for $\varphi$.*

Before proving the proposition we first need a few technical results.

LEMMA 3.8. *Let $N \subset X$ a compact subset, and suppose $\{x_n\} \subset N$ is a sequence such that $x_n \to x \in N$ and $\{\tau_n\} \subset [0, \infty)$ is a sequence such that $\tau_n \to \infty$ as $n \to \infty$. If there exist backward orbit segments $\gamma_{x_n}^- : [-\tau_n, 0] \to N$ through $x_n$ for each $n \in \mathbb{N}$, then there exists a complete backward orbit $\gamma_x^- : \mathbb{T}^- \to N$ through $x$, such that $\gamma_{x_n}^-|_I \to \gamma_x^-|_I$ uniformly on every compact interval $I \subset \mathbb{T}^-$.*

PROOF. Let $x_n^0 = x_n$ and $x^0 = x$. Choose $n_1$ such that $\tau_n > 1$ for $n \geq n_1$ and define $x_n^1 = \gamma_{x_n}^-(-1)$. By choosing a subsequence if necessary, denoted again by $x_n^1$, we can assume by the compactness of $N$ that $x_n^1 \to x^1$. There exists a unique forward orbit $\gamma_{x^1}^+$ through $x^1$, and the function $\gamma_{x^0}^-$ defined by $\gamma_{x^0}^-(t) = \gamma_{x^1}^+(1+t)$ on $[-1,0]$ is a piece of a backward orbit through $x^0$. We can repeat this process to recursively obtain integers $n_k$ such that $\tau_n > k$ for $n \geq n_k$ and sequences $x_n^k \to x^k$. Now glue together the backward orbits $\gamma_{x^k}^-$ defined on $[-1,0]$ through $x^k$ to obtain a backward orbit $\gamma_x^-$ through $x$ defined on $(-\infty,0]$. Without loss of generality fix $\tau > 0$ and consider the interval $I = [-\tau, 0]$. Let $y_n = \gamma_{x_n}^-(-\tau)$. For $n$ sufficiently large the orbit segments $\gamma_{x_n}^-|_I$ are defined, and are equivalent to the forward orbit segments $\gamma_{y_n}^+ : [0,\tau] \to N$ under time reversal. The convergence of the forward orbit segments is due to the uniform continuity of $\varphi$ on compact sets, from which the convergence of the backward orbits follows. ∎

LEMMA 3.9. *Suppose $N$ is a compact set with $S = \mathrm{Inv}(N) \subset \mathrm{int}\,(N)$. For $\tau > 0$ define $N_\tau^+ = \{x \in N \mid \varphi([0,\tau],x) \subset N\}$. Then there exists $\delta_\tau > 0$ such that $B_{\delta_\tau}(S) \subset N_\tau^+$, i.e. $\varphi([0,\tau], B_{\delta_\tau}(S)) \subset N$.*

PROOF. By continuity and the invariance of $S = \mathrm{Inv}(N)$, for each $x \in S$ there exists $\delta_x > 0$ such that $\varphi([0,\tau], B_{\delta_x}(x)) \subset N$. By compactness there exists finitely many such balls $B_{\delta_{x_i}}(x_i)$ such that $S \subset \cup_i B_{\delta_{x_i}}(x_i)$. Since this union is open, there exists $\delta > 0$ such that $B_\delta(S) \subset \cup_i B_{\delta_{x_i}}(x_i) \subset N_\tau^+$. ∎

LEMMA 3.10. *Suppose $N$ is a compact set such that $\mathrm{Inv}(N) \subset \mathrm{int}\,(N)$ with the property that there are no backward orbits $\gamma_x^- : \mathbb{T}^- \to N$ through $x \in N \setminus \mathrm{Inv}(N)$. Then for every $\epsilon > 0$ there exists $\tau > 0$ such that there are no backward orbit segments $\gamma_x^- : [-\tau, 0] \to N$ through $x \in N \setminus B_\epsilon(\mathrm{Inv}(N))$.*

PROOF. Define $N_\tau^- = \{x \in N \mid \exists$ a backward orbit segment $\gamma_x^- : [-\tau,0] \to N\}$ for $\tau > 0$, and let $\epsilon_\tau^- = \sup_{x \in N_\tau^-}\{d(x, \mathrm{Inv}(N))\}$. First we show that $\epsilon_\tau^- \to 0$ as $\tau \to \infty$. Suppose not, then there exist sequences $\tau_n \to \infty$ and $x_n \in N$ with $x_n \to x \in N \setminus \mathrm{Inv}(N)$ and $\gamma_{x_n}^- : [-\tau_n, 0] \to N$. By Lemma 3.8, there exists a backward orbit $\gamma_x^- \subset N$, which is a contradiction. Finally choose $\tau > 0$ large enough so that $\epsilon_\tau^- < \epsilon$. ∎

LEMMA 3.11. *An invariant set $A$ is an attractor if and only if there exists a compact neighborhood $N$ of $A$ such that there are no backward orbits $\gamma_x^- : \mathbb{T}^- \to N$ through $x \in N \setminus A$. In particular, $A = \mathrm{Inv}(N)$.*

PROOF. Suppose $A$ is an attractor. A compact trapping region $N$ such that $A = \omega(N)$ exists by Lemma 3.17. If there is a backward orbit $\gamma_x^- \subset N$ with $x \in N \setminus A$, then by Proposition 2.11(vii), $x \in \omega(N) = A$, a contradiction. The compact trapping region $N$ is a neighborhood with the required property.

Now suppose $N$ is a compact neighborhood of $A$ with the the property that no backward orbits $\gamma_x^- : \mathbb{T}^- \to N$ exist through $x \in N \setminus A$. The latter and the fact that $N$ is a compact neighborhood of $A$, implies that $A = \mathrm{Inv}(N) \subset \mathrm{int}\,(N)$. By Lemma 3.9, there exists $\delta_1 > 0$ such that $B_{\delta_1}(A) \subset N$ and $\varphi([0,1], B_{\delta_1}(A)) \subset N$. Fix $\epsilon < \delta_1$ and consider the neighborhood $B_\epsilon(A)$. By Lemma 3.10, there exists $\tau > 0$ such that there are no backward orbit segments $\gamma_x^- : [-\tau, 0] \to N$ for $x \in N \setminus B_\epsilon(A)$.

Now we construct an attracting neighborhood for $A$. By Lemma 3.9, there exists $\delta_{2\tau} > 0$ such that the neighborhood $U = B_{\delta_{2\tau}}(A)$ satisfies $\varphi([0, 2\tau], U) \subset N$. We show that $U$ is an attracting neighborhood. Let $V_0 = \varphi(\tau, U) \subset N$. For each $x \in V_0$ there exists a backward orbit segment $\gamma_x^- : [-\tau, 0] \to N$. The definition of $\tau$ from Lemma 3.10 implies that $V_0 \subset B_\epsilon(A)$. Moreover, by our choice of $\epsilon < \delta_1$ we have $\varphi([0,1], B_\epsilon(A)) \subset N$ so that $V_1 = \varphi([0,1], V_0) \subset N$. Thus for each $x \in V_1$ there exists a backward orbit segment $\gamma_x^- : [-\tau - 1, 0] \to N$ which implies that $\gamma_x^-([-\tau, 0]) \subset N$ so that $V_1 \subset B_\epsilon(A)$. We can repeat this argument inductively to prove that $V_k = \varphi([k-1, k], V_0) \subset B_\epsilon(A) \subset N$ for all $k > 0$, and thus $\varphi([0,\infty), U) \subset N$. Finally, $A \subset U$ implies that $A = \omega(A) \subset \omega(U)$, and $\varphi([0,\infty), U) \subset N$ implies that $\omega(U) \subset \mathrm{Inv}(N) = A$. Therefore $A = \omega(U)$, and $A$ is an attractor by Proposition 3.5. ■

*Proof of Proposition 3.7.* Let $N$ be a (compact) attracting neighborhood for $A$ in $X$. There exists a neighborhood $N' \subset N$ of $A'$ in $X$ such that $N' \cap A$ is an attracting neighborhood for $A'$ in $A$. Indeed, since $A'$ an attractor in $A$ there exists a neighborhood $\widetilde{N}' \supset A'$ in $A$, such that $\omega(\widetilde{N}', \varphi|_A) = A'$. Choose $N' \subset N$ such that $A' \subset N' \cap A \subset \widetilde{N}'$, then

$$A' = \omega(A', \varphi|_A) \subset \omega(N' \cap A, \varphi|_A) \subset \omega(\widetilde{N}', \varphi|_A) = A',$$

which shows that $N' \cap A$ is an attracting neighborhood for $A'$ in $A$.

Suppose that $\gamma_x^-$ is a backward orbit through $x \in N' \setminus A'$ such that $\gamma_x^- \subset N'$. By Proposition 2.11(vii) this implies that $x \in \omega(N') \subset \omega(N) = A$. Indeed the same holds for all $y \in \gamma_x^-$, which yields $\gamma_x^- \subset A$, and therefore $\gamma_x^- \subset N' \cap A$. Proposition 2.11(vii) again implies that $x \in \omega(N' \cap A) = A'$, a contradiction. The criterion in Lemma 3.11 the reveals that $A'$ is an attractor for $\varphi$ in $X$. ■

### 3.3. Repellers.

A *repelling region* is a backward invariant set $U \subset X$ such that

$$\varphi(\tau, \mathrm{cl}\,(U)) \subset \mathrm{int}\,(U) \quad \text{for some } \tau < 0.$$

A set $R \subset X$ is called an *repeller* if there exists a repelling region $U$ such that $R = \mathrm{Inv}^+(U)$.

As attractors are characterized by attracting neighborhoods, a similar characterization can be given for repellers.

PROPOSITION 3.12. *A repeller $R$ is a compact, forward-backward invariant set. If $U$ is a repelling region for $R$, then $R = \mathrm{Inv}^+(U) = \mathrm{Inv}^+(\mathrm{cl}\,(U)) = \alpha(U)$ and $R \subset \mathrm{int}\,(U)$.*

PROOF. By Lemma 2.10, repellers are forward-backward invariant. Let $U$ be a repelling region for $R$. If $R' = \mathrm{Inv}^+(\mathrm{cl}\,(U))$, then $\varphi(t, R') \subset R'$, for all $t \geq 0$, and in particular $\varphi(-\tau, R') \subset R'$ ($\tau < 0$ for the definition of repelling region). Therefore, $R' \subset \varphi(\tau, \varphi(-\tau, R')) \subset \varphi(\tau, R') \subset \varphi(\tau, \mathrm{cl}\,(U)) \subset \mathrm{int}\,(U) \subset U$. This implies $R' \subset \mathrm{Inv}^+(U) = R \subset \mathrm{Inv}^+\big(\mathrm{cl}\,(U)\big) = R'$, which shows that $R' = R$.

Since $R \in \mathsf{Invset}^{\pm}(\varphi)$, Proposition 2.13(viii) implies that $\mathrm{cl}\,(R) \subset \alpha(R) = R \subset \mathrm{cl}\,(R)$, and thus $\mathrm{cl}\,(R) = R$. Furthermore, $R = \alpha(R) \subset \alpha(U) \subset \alpha\big(\mathrm{cl}\,(U)\big) = \mathrm{Inv}^+\big(\mathrm{cl}\,(U)\big) = R$, which proves that $R = \alpha(U)$. ■

While the intersection of invariant sets need not be invariant, Lemma 2.9 immediately implies the following result.

PROPOSITION 3.13. *If $A \in \mathsf{Att}(\varphi)$ and $R \in \mathsf{Rep}(\varphi)$, then $A \cap R \in \mathsf{Invset}(\varphi)$.*

A neighborhood $U \subset X$ is a *repelling* neighborhood if

$$\alpha(U) \subset \mathrm{int}\,(U). \tag{5}$$

The set of repelling neighborhoods in $X$ is denoted by $\mathsf{RNbhd}(\varphi)$. Via repelling neighborhoods the following characterization of repellers holds.

PROPOSITION 3.14. *A set $R \subset X$ is a repeller if and only if there exists a neighborhood $U$ of $R$ such that $R = \alpha(U)$. In particular, $U$ is a repelling neighborhood. Moreover, for every repelling neighborhood $U$ there exists a repelling region $U' \subset U$.*

COROLLARY 3.15. *Let $U \subset X$ be a repelling neighborhood, then $R = \mathrm{Inv}^+(U) = \mathrm{Inv}^+(\mathrm{cl}\,(U)) = \alpha(U)$ is a repeller.*

We postpone the proofs of Proposition 3.14 and Corollary 3.15 and a result for repellers analogous to Proposition 3.7 to Section 3.5, since the arguments are greatly simplified once the duality between attractors and repellers is established.

**3.4. Attractor-repeller pairs.** For an attractor $A \in \mathsf{Att}(\varphi)$, with trapping region $U$, the *dual repeller* of $A$ is defined by

$$A^* = \mathrm{Inv}^+(U^c),$$

For a repeller $R \in \mathsf{Rep}(\varphi)$, with repelling region $U$, the *dual attractor* of $R$ is defined by

$$R^* = \mathrm{Inv}(U^c).$$

PROPOSITION 3.16. *For an attractor $A$, the dual repeller $A^*$ is compact, forward-backward invariant, and characterized by*

$$A^* = A^{\oplus} = \{x \in X \mid \omega(x) \cap A = \varnothing\}, \tag{6}$$

*and therefore independent of the chosen trapping region for $A$. Moreover, if $X$ is invariant, then $A^*$ is strongly invariant. For a repeller $R$, the dual attractor $R^*$ is compact, invariant, and characterized by*

$$R^* = R^{\ominus} = \left\{x \in X \;\middle|\; \exists\, \gamma_x \text{ such that } \alpha_{\mathrm{o}}(\gamma_x^-) \cap R = \varnothing \right\}, \tag{7}$$

*and therefore independent of the chosen repelling region for $R$.*

PROOF. Let $U$ be trapping region for $A$, and let $A^* = \mathrm{Inv}^+(U^c)$. Since $A^* \subset U^c$, we have $\mathrm{cl}\,(A^*) \subset \mathrm{cl}\,(U^c)$, and $\mathrm{cl}\,(A^*)$ is forward invariant by Remark 2.12. We show that $\mathrm{cl}\,(A^*) \subset U^c$. Indeed, if $x \in \mathrm{cl}\,(A^*) \cap U$, then since $U$ is a trapping region, $\varphi(\tau, x) \in \mathrm{int}\,(U)$ for some $\tau > 0$. However, $\varphi(\tau, x) \in \mathrm{cl}\,(A^*) \subset \mathrm{cl}\,(U^c)$, which is a contradiction. Therefore, $A^* \subset \mathrm{cl}\,(A^*) \subset \mathrm{Inv}^+(U^c) = A^*$, so that $A^* = \mathrm{cl}\,(A^*)$, and hence $A^*$ is compact.

By Proposition 2.16 the set $A^{\oplus}$ is forward-backward invariant. If $x \in U$, then, since $U$ is a trapping region, $\omega(x) \subset \omega(U) = A$, and therefore $A^{\oplus} \subset U^c$. Since $A^{\oplus}$ is forward-backward invariant, $A^{\oplus} \subset \mathrm{Inv}^+(U^c) = A^*$. If $x \in A^*$, then $\omega(x) \subset \omega(A^*)$ and by Proposition 2.11(iii), $\omega(A^*) = \mathrm{Inv}(A^*) \subset A^*$, so that $\omega(x) \subset A^* \subset U^c$. Therefore $\omega(x) \cap A = \varnothing$, and thus $x \in A^{\oplus}$. Hence $A^* \subset A^{\oplus}$. Combining these inclusions, $A^* = A^{\oplus}$.

Let $U$ be a repelling region for $R$ and let $R^* = \mathrm{Inv}(U^c)$. By definition $R^* \subset U^c$, and thus $\mathrm{cl}\,(R^*) \subset \mathrm{cl}\,(U^c)$. Since $R^*$ is invariant by definition, $\mathrm{cl}\,(R^*)$ is compact and invariant, see Remark 2.12. If $x \in \mathrm{cl}\,(R^*) \cap U$, then, since $U$ is a repelling region, $\varphi(-\tau, x) \subset \mathrm{int}\,(U)$ for some $\tau > 0$, i.e. for all backward orbits $\gamma_x^-$ it holds that $\gamma_x^-(-\tau) \in \mathrm{int}\,(U)$. However, since $\mathrm{cl}\,(R^*)$ is invariant, there exists a backward orbit $\gamma_x^- \subset \mathrm{cl}\,(R^*) \subset \mathrm{cl}\,(U^c)$, which is a contradiction. This shows that $\mathrm{cl}\,(R^*) \subset U^c$, and therefore $\mathrm{cl}\,(R^*) \subset \mathrm{Inv}(U^c) = R^*$, which proves that $R^* = \mathrm{cl}\,(R^*)$ and $R^*$ is compact.

By Proposition 2.16 the set $R^{\ominus}$ is invariant. If $x \in U$, then there are two possibilities. Either there is no complete orbit $\gamma_x$ through $x$ in $X$, in which case $x \notin R^{\ominus}$, or since $U$ is a repelling region, every complete orbit $\gamma_x$ has the property that $\alpha_{\mathrm{o}}(\gamma_x^-) \subset U$. The latter statement follows from the fact that $\alpha_{\mathrm{o}}(\gamma_x^-)$ is invariant by Proposition 2.15. Indeed, by definition $\alpha_{\mathrm{o}}(\gamma_x^-) \subset \mathrm{cl}\,(U)$; since $U$ is a repelling region, there exists $\tau > 0$ such that $\varphi(-\tau, \mathrm{cl}\,(U)) \subset \mathrm{int}\,(U) \subset U$, and hence $\alpha_{\mathrm{o}}(\gamma_x^-) \subset \varphi(-\tau, \varphi(\tau, \alpha_{\mathrm{o}}(\gamma_x^-))) = \varphi(-\tau, \alpha_{\mathrm{o}}(\gamma_x^-)) \subset U$. Moreover, $\alpha_{\mathrm{o}}(\gamma_x^-) \subset R$, and thus $x \notin R^{\ominus}$. Therefore $R^{\ominus} \subset U^c$. Consequently, since $R^{\ominus}$ is invariant, $R^{\ominus} \subset \mathrm{Inv}(U^c) = R^*$. If $x \in R^*$, then by invariance there exists an orbit $\gamma_x \subset R^*$. By compactness, $\alpha_{\mathrm{o}}(\gamma_x^-) \subset R^* \subset U^c$, which implies that $x \in R^{\ominus}$, and thus $R^* \subset R^{\ominus}$. Combining these inclusions gives that $R^* = R^{\ominus}$. ∎

This proposition allows us to make the following fundamental definition. Given $A \in \mathsf{Att}(\varphi)$ or $R \in \mathsf{Rep}(\varphi)$, the pair $(A, A^*)$, or equivalently $(R^*, R)$, is called an *attractor-repeller pair for $\varphi$.*

LEMMA 3.17. *If $U$ is a trapping region, then $N = \mathrm{cl}\,(U)$ is also a trapping region, and $N^c$ is a repelling region. If $V$ is a repelling region, then $\mathrm{int}\,(V)$ is also a repelling region, and $V^c$ is a trapping region.*

PROOF. Suppose $U$ is a trapping region and define $N = \mathrm{cl}\,(U)$. Since there exists $\tau > 0$ such that $\varphi(\tau, \mathrm{cl}\,(U)) \subset \mathrm{int}\,(U)$, we have $\varphi(\tau, \mathrm{cl}\,(N)) = \varphi(\tau, \mathrm{cl}\,(U)) \subset \mathrm{int}\,(U) = \mathrm{int}\,(N)$. By Remark 2.12 the closure of a forward invariant set is also forward invariant, so $N$ is a trapping region.

By Proposition 2.5, $N^c = (\mathrm{cl}\,(U))^c = \mathrm{int}\,(U^c)$ is backward invariant. Since $N$ is a trapping region, there exists a $\tau > 0$ such that $\varphi(\tau, N) \subset \mathrm{int}\,(N)$. Therefore, $N \subset \varphi(-\tau, \varphi(\tau, N)) \subset \varphi(-\tau, \mathrm{int}\,(N))$ and by taking complements

$$N^c \supset \varphi(-\tau, \mathrm{int}\,(N))^c = \varphi\big(-\tau, (\mathrm{int}\,(N))^c\big) = \varphi(-\tau, \mathrm{cl}\,(N^c)).$$

Since $N^c$ is open, $\mathrm{int}\,(N^c) = N^c$, which gives $\varphi(-\tau, \mathrm{cl}\,(N^c)) \subset N^c = \mathrm{int}\,(N^c)$, and thus $N^c$ is a repelling region.

Let $V$ be a repelling region, then, by Proposition 2.5, $V^c$ is forward invariant. By assumption $\varphi(\tau, \mathrm{cl}\,(V)) \subset \mathrm{int}\,(V)$, for some $\tau < 0$ and therefore $(\mathrm{int}\,(V))^c \subset \varphi(\tau, \mathrm{cl}\,(V))^c = \varphi(\tau, (\mathrm{cl}\,(V))^c)$ and

$$\varphi\big(-\tau, (\mathrm{int}\,(V))^c\big) \quad \subset \varphi(-\tau, \varphi(\tau, (\mathrm{cl}\,(V))^c)) \subset (\mathrm{cl}\,(V))^c = \mathrm{int}\,(V^c).$$

Since $(\mathrm{int}\,(V))^c = \mathrm{cl}\,(V^c)$, we obtain $\varphi(-\tau, \mathrm{cl}\,(V^c)) \subset \mathrm{int}\,(V^c)$, which proves that $V^c$ is a trapping region.

Let $W = \mathrm{int}\,(V)$. Since there exists $\tau < 0$ such that $\varphi(\tau, \mathrm{cl}\,(V)) \subset \mathrm{int}\,(V)$, we have $\varphi(\tau, \mathrm{cl}\,(W)) \subset \varphi(\tau, \mathrm{cl}\,(V)) \subset \mathrm{int}\,(V) = \mathrm{int}\,(W)$. Moreover, $V^c$ is forward invariant by Proposition 2.5, and hence $\mathrm{cl}\,(V^c)$ is forward invariant by Remark 2.12. Thus $W = (\mathrm{cl}\,(V^c))^c$ is backward invariant by Proposition 2.5. Therefore, $W$ is a repelling region. ■

PROPOSITION 3.18. *For an attractor $A$, the dual repeller $A^*$ is a repeller, and for a repeller $R$, the dual attractor $R^*$ is an attractor. Moreover $(A^*)^* = A$ and $(R^*)^* = R$.*

PROOF. Let $U$ be a trapping for $A$, then by Proposition 3.1, $A = \mathrm{Inv}(U) = \mathrm{Inv}(\mathrm{cl}\,(U))$. By Lemma 3.17, $N = \mathrm{cl}\,(U)$ is also a trapping region, and since $A^*$ is independent of the chosen trapping region by Proposition 3.16, we conclude that $A^* = \mathrm{Inv}^+(U^c) = \mathrm{Inv}^+(N^c)$. By Lemma 3.17, $N^c$ is a repelling region and therefore $A^*$ is a repeller. Since $N^c$ is a repelling region and $(N^c)^c = N$ is a trapping region, $(A^*)^* = \mathrm{Inv}((N^c)^c) = \mathrm{Inv}(N) = A$.

Let $U$ be a repelling region with $R = \mathrm{Inv}^+(U)$. By definition, $R^* = \mathrm{Inv}(U^c)$. By Lemma 3.17, $U^c$ is a trapping region. Therefore, $R^*$ is an attractor by definition. Since $U^c$ is a trapping region and $(U^c)^c = U$ is repelling region, $(R^*)^* = \mathrm{Inv}^+((U^c)^c) = \mathrm{Inv}^+(U) = R$. ■

The follow result gives an important characterization of attractor-repeller pairs which characterizes the dynamics by a simple partial order on the attractor-repeller pair.

THEOREM 3.19. *For $A, R \subset X$ the following statements are equivalent.*

(i) *$(A, R)$ is an attractor-repeller pair.*
(ii) *$A$ and $R$ are disjoint, compact sets with $A \in \mathsf{Invset}(\varphi)$ and $R \in \mathsf{Invset}^+(\varphi)$ such that for every $x \in X \setminus (A \cup R)$ and every backward orbit $\gamma_x^-$ through $x$ we have $\alpha_{\rm o}(\gamma_x^-) \subset R$ and $\omega(x) \subset A$.*

PROOF. Suppose $(A, R)$ is an attractor-repeller pair, then $R = A^*$ and $A \in \mathsf{Invset}(\varphi)$, $R \in \mathsf{Invset}^\pm(\varphi)$ and both sets are compact by Proposition 3.1. By the definition of dual repeller it follows that $A \cap A^* = \varnothing$. By Lemma 3.4 there exists $\epsilon_0 > 0$ such that $\omega(B_\epsilon(A)) = A$ for all $0 < \epsilon < \epsilon_0$. By Equation (6) if $x \notin A \cup R$, then $\omega(x) \cap A \neq \varnothing$, which implies that there exists $\tau > 0$ such that $d(\varphi(\tau, x), A) < \epsilon_0$ by the definition of omega-limit set. Hence $\omega(x) = \omega(\varphi(\tau, x)) \subset \omega(B_{\epsilon_0}(A)) = A$. By Proposition 3.5 and Lemma 3.17 there exists a compact trapping region $N$ with $A = \mathrm{Inv}(N)$ and $N \subset B_{\min\{\epsilon_0, d(x,A)\}}(A)$. Lemma 3.17 yields that $N^c$ is a repelling region and $R = \mathrm{Inv}^+(N^c)$. . Hence $x \in N^c$ implies $\gamma_x^- \subset N^c$ whenever such a backward orbit exists. Moreover, since $\alpha_{\rm o}(\gamma_x^-)$ is invariant with $\alpha_{\rm o}(\gamma_x^-) \subset \mathrm{cl}\,(N^c)$ and $N^c$ is a repelling region, we have $\alpha_{\rm o}(\gamma_x^-) \subset \mathrm{Inv}(N^c) \subset \mathrm{Inv}^+(N^c) = R$, which proves that (i) implies (ii).

Suppose $A$ and $R$ are disjoint, compact, forward invariant sets in $X$, and $\omega(x) \subset A$ for every $x \in X \setminus (A \cup R)$, then $R = A^\oplus$. Indeed, if $x \in R$, then by forward invariance and compactness of $R$, we have $\omega(x) \subset \omega(R) \subset R$. Since $A \cap R = \varnothing$ this implies $\omega(x) \cap A = \varnothing$, and therefore $R \subset A^\oplus$. Conversely, if $x \in A^\oplus$, then $\omega(x) \cap A = \varnothing$. By assumption this implies that $x \in A \cup R$. Moreover, $A$ is compact and forward invariant and thus $x \in A$, implies $\omega(x) \subset \omega(A) \subset A$. Consequently, $x \in A^\oplus$ implies $x \in R$ and therefore $A^\oplus \subset R$. Summarizing $R = A^\oplus$.

If $A$ is an attractor, then Equation (6) implies $R = A^*$ is its dual repeller, and hence $(A, R)$ is an attractor-repeller pair. We now show that $A$ is an attractor. Fix a compact neighborhood $U$ of $A$ such that $A \subset \mathrm{int}\,(U)$ and $U \cap R = \varnothing$. If $x \in U \setminus A$ and $\gamma_x^-$ is a backward orbit through $x$, then by (ii), since $x \in X \setminus (A \cup R)$, we have $\alpha_o(\gamma_x^-) \subset R$, and hence $\gamma_x^- \not\subset U$. Lemma 3.11 implies that $A$ is an attractor. ∎

REMARK 3.20. Note that in Condition (ii) of Theorem 3.19 $R$ is not a priori required to be forward-backward invariant. The combination with the other assumptions in (ii) implies that $R = A^\oplus$, and therefore $R$ is forward-backward invariant as a consequence. If we assume in (ii) that both $A$ and $R$ are forward invariant, then $(\mathrm{Inv}(A), R)$ is an attractor-repeller pair.

**3.5. Consequences of duality.** We can use the characterization of attractor-repeller pairs in Theorem 3.19 to give an alternative characterization of attracting neighborhoods of an attractor in terms of the dual repeller.

PROPOSITION 3.21. *Let $A \in \mathsf{Att}(\varphi)$. A subset $U \subset X$ is an attracting neighborhood for $A$ if and only if $U$ is a neighborhood of $A$ and* $\mathrm{cl}\,(U) \cap A^* = \varnothing$.

PROOF. Let $U$ be an attracting neighborhood with $A = \omega(U)$. Then by Proposition 3.5, we have $A \subset \mathrm{int}\,(U)$. By Proposition 2.11(vi), we have $A = \omega(\mathrm{cl}\,(U))$ and by Lemma 3.2, $\mathrm{cl}\,(U)$ is an attracting neighborhood for $A$. Therefore, $\omega(x) \subset A$ for all $x \in \mathrm{cl}\,(U)$, and $\mathrm{cl}\,(U) \cap A^* = \varnothing$ by Equation (6).

Let $A \subset \mathrm{int}\,(U)$ and $\mathrm{cl}\,(U) \cap A^* = \varnothing$. By Theorem 3.19, we have $\omega(x) \subset A$ for all $x \in \mathrm{cl}\,(U)$. Thus by continuity, for each $x \in \mathrm{cl}\,(U)$ there exists $\tau_x \geq 0$ and $\delta_x > 0$ such that $\varphi(t, B_{\delta_x}(x)) \subset U$ for all $t \geq \tau_x$. By compactness, there exists $\tau \geq 0$ such that $\varphi(t, x) \in U$ for all $t \geq \tau$ and $x \in \mathrm{cl}\,(U)$, i.e. $\varphi(t, \mathrm{cl}\,(U)) \subset U$ for all $t \geq \tau$. By Proposition 2.11(iii) and (vi) this implies that $\omega(U) = \omega(\mathrm{cl}\,(U)) = \mathrm{Inv}(\mathrm{cl}\,(U)) \subset \mathrm{cl}\,(U)$. In particular $\omega(U) \cap A^* = \varnothing$. Now $(\omega(U), A^*)$ is a pair of disjoint, compact sets with $\omega(U) \in \mathsf{Invset}(\varphi)$ and $A^* \in \mathsf{Invset}^{\pm}(\varphi)$ which satisfy property (ii) in Theorem 3.19. Indeed, $A \subset \mathrm{int}\,(U) \subset U$, which implies $A = \omega(A) \subset \omega(U)$, and since $(A, A^*)$ satisfies (ii), also $(\omega(U), A^*)$ satisfies this property. Theorem 3.19 implies that $(\omega(U), A^*)$ is an attractor-repeller pair, and hence $A = \omega(U)$ by Proposition 3.18. Therefore, $U$ is an attracting neighborhood for $A$, since $\omega(U) = A \subset \mathrm{int}\,(U)$. ■

COROLLARY 3.22. *The map* $\mathrm{Inv}\colon \mathsf{ANbhd}(\varphi) \to \mathsf{Att}(\varphi)$, *restricted to the compact attracting neighborhoods, is continuous in the Hausdorff metric. Indeed, it is locally constant.*

PROOF. Let $U, V$ be compact sets in $\mathsf{ANbhd}(\varphi)$ and define $d(U, V) = \mathrm{dist}_{\,\mathrm{Hausdorff}}(U, V)$. Let $A = \mathrm{Inv}(U)$ be the attractor in $U$. By Proposition 3.21, $U \cap A^* = \varnothing$. Since $U$ and $A^*$ are compact, $\mathrm{dist}\,(U, A^*) = \delta_1 > 0$. Also, since $U$ is a neighborhood of $A$, there exists $\delta_2 > 0$ such that $B_{\delta_2}(A) \subset U$. Let $\delta = \min\{\delta_1/2, \delta_2/2\}$. Then for any neighborhood $V$ with $d(U, V) < \delta$ we have that $A \subset V$ and $V \cap A^* = \varnothing$. Therefore, $A = \mathrm{Inv}(V)$ by Proposition 3.21, which proves the result. ■

By a similar argument, $\mathrm{Inv}^+ : \mathsf{RNbhd}(\varphi) \to \mathsf{Rep}(\varphi)$ is also locally constant on the compact repelling neighborhoods, using Proposition 3.25 below.

LEMMA 3.23. *If $U \subset X$ is an attracting neighborhood for an attractor $A = \omega(U)$, then $U^c$ is a neighborhood of $A^*$, with $\mathrm{cl}\,(U^c) \cap A = \varnothing$ and $A^* = \mathrm{Inv}^+(U^c) = \alpha(U^c)$.*

PROOF. By Proposition 3.21 an attracting neighborhood $U$ of $A$ satisfies $\mathrm{cl}\,(U) \cap A^* = \varnothing$. It follows that $A^* \subset (\mathrm{cl}\,(U))^c$ and since $\mathrm{int}\,(U^c) = (\mathrm{cl}\,(U))^c$, we have $A^* \subset \mathrm{int}\,(U^c) \subset U^c$. This shows that $U^c$ is a neighborhood of $A^*$. Since $A \subset \mathrm{int}\,(U)$, it follows that $A \cap (\mathrm{int}\,(U))^c = \varnothing$, and by $\mathrm{cl}\,(U^c) = (\mathrm{int}\,(U))^c$ we derive that $\mathrm{cl}\,(U^c) \cap A = \varnothing$. Let $U' \subset U$ be a trapping region for $A$. By Lemma 3.17, $U'$ can taken to be a compact trapping region and $U'^c$ is a repelling region. By definition, $A^* = \mathrm{Inv}^+(U'^c)$ and by Proposition 3.12 $A^* = \mathrm{Inv}^+(U'^c) = \alpha(U'^c) \supset \alpha(U^c)$. Since $A^* \subset U^c$, we have by Proposition 2.13(viii) that $A^* = \alpha(A^*) \subset \alpha(U^c)$, which proves that $A^* = \alpha(U^c)$. Since $\alpha(U^c) \subset U^c$, Proposition 2.13(vi) implies $\alpha(U^c) = \mathrm{Inv}^+(U^c)$. ■

Now recall Proposition 3.14 which states that a set $R \subset X$ is a repeller if and only if there exists a neighborhood $U$ of $R$ such that $R = \alpha(U)$. In particular, $U$ is a repelling neighborhood and for every $U$ there exists a repelling region $U' \subset U$.

*Proof of Proposition 3.14.* If $R$ is a repeller, then there exists a repelling region $U$ such that $R = \mathrm{Inv}^+(U) = \alpha(U) \subset \mathrm{int}\,(U)$ by Proposition 3.12.

Let $U \subset X$ be a repelling neighborhood such that $R = \alpha(U) \subset \mathrm{int}\,(U)$. We now show that $R$ is a repeller. By Proposition 2.13(i), $R$ is compact and $R \in \mathsf{Invset}^+(\varphi)$. Since $R \subset \mathrm{int}\,(U)$, it follows that $R^c \supset \big(\mathrm{int}\,(U)\big)^c = \mathrm{cl}\,(U^c)$. This implies that we have two compact sets $\mathrm{cl}\,(U^c)$ and $R$, such that

$$\mathrm{cl}\,(U^c) \cap R = \varnothing.$$

Since $X$ is a compact metric space, there exists a compact neighborhood $V$ for $R$, such that $\mathrm{cl}\,(U^c) \cap V = \varnothing$. The latter implies that $V \subset \big(\mathrm{cl}\,(U^c)\big)^c = \mathrm{int}\,(U)$, and hence

$$R \subset \mathrm{int}\,(V) \subset V \subset \mathrm{int}\,(U) \subset U. \tag{8}$$

Since $R$ is compact and $R \in \mathsf{Invset}^+(\varphi)$, Proposition 2.13(iv) and (viii) imply that

$$R = \mathrm{cl}\,(R) \subset \alpha(R) \subset \alpha(V) \subset \alpha(U) = R,$$

and thus $R = \alpha(V) \subset \mathrm{int}\,(V)$. By the definition of alpha-limit set there exists a $\tau < 0$ such that $\mathrm{cl}\,(\varphi((-\infty, \tau], V)) \subset \mathrm{int}\,(V)$, and thus, since $V$ is compact,

$$\varphi(t, \mathrm{cl}\,(V)) = \varphi(t, V) \subset \mathrm{cl}\,(\varphi(t, V)) \subset \mathrm{int}\,(V), \quad \forall\, t \leq \tau.$$

The latter yields $(\mathrm{int}\,(V))^c \subset \varphi(t, \mathrm{cl}\,(V))^c = \varphi(t, (\mathrm{cl}\,(V))^c)$, for $t \leq \tau < 0$, and consequently,

$$\varphi\big(-t, (\mathrm{int}\,(V))^c\big) \subset \varphi\big(-t, \varphi(t, (\mathrm{cl}\,(V))^c)\big) \subset (\mathrm{cl}\,(V))^c = \mathrm{int}\,(V^c).$$

Since $(\mathrm{int}\,(V))^c = \mathrm{cl}\,(V^c)$, we obtain $\varphi(-t, \mathrm{cl}\,(V^c)) \subset \mathrm{int}\,(V^c)$, for all $-t \geq -\tau > 0$. Lemma 3.2 yields that $V^c$ is an attracting neighborhood, with $\omega(V^c) \subset \mathrm{int}\,(V^c)$. By Proposition 3.5, $A = \omega(V^c)$ is an attractor, and by Lemma 3.23, $A^* = \alpha(V) = R$, which proves that $R$ is a repeller.

Finally we show that there exists a repelling region $U' \subset U$, for every repelling neighborhood $U$ for $R$. From (8) we obtain

$$U^c \subset \mathrm{cl}\,(U^c) \subset V^c \subset \mathrm{cl}\,(V^c) \subset R^c,$$

and we define $W = \mathrm{cl}\,(V^c)$, which is a compact neighborhood of $A = \omega(V^c) = \omega(W)$ and $W \cap R = \varnothing$. Proposition 3.21 yields that $W$ is an attracting neighborhood for $A$, and by Lemma 3.2 there exists a $\tau > 0$, such that $\varphi(t, W) \subset \mathrm{int}\,(W)$ for all $t \geq \tau$. By construction $W' = \varphi([0, \tau], W) \supset W$ is a compact, forward invariant set, cf. Proposition 3.5. Another property of $W'$ is that $W' \cap R = \varnothing$. Indeed, if not, then there would exist a point $x \in W$ and a time $0 < t_0 < \tau$, such that $\varphi(t_0, x) \in R$. Furthermore, $\varphi(t, \varphi(t_0, x)) = \varphi(t + t_0, x) \in \mathrm{int}\,(W)$ for all $t \geq \tau - t_0$. However, since $R$ is forward invariant, $\varphi(t, \varphi(t_0, x)) \in R$ for all $t \geq 0$, which is a contradiction. By Proposition 3.21 we conclude that $W'$ is a attracting neighborhood. Combining this with the fact that $W'$ is forward invariant, Lemma 3.3 yields that $W'$ is a trapping region for $A$. By Lemma 3.17 and Lemma 3.23 we have that $U' = (W')^c$ is a repelling region for $R = A^*$. Since $U^c \subset W \subset W'$, it follows that $U' \subset U$, which completes the proof. ■

COROLLARY 3.24. *A set $U$ is an attracting neighborhood if and only if $U^c$ is a repelling neighborhood.*

PROOF. Combine Proposition 3.14 and Lemma 3.23. ■

PROPOSITION 3.25. *Let $R \in \mathsf{Rep}(\varphi)$ be a repeller. A subset $U \subset X$ is an repelling neighborhood for $R$ if and only if $U$ is a neighborhood of $R$ and $\mathrm{cl}\,(U) \cap R^* = \varnothing$.*

PROOF. If $U$ is a repelling neighborhood for $A^* = R = \alpha(U) \subset \mathrm{int}\,(U)$, then $\mathrm{cl}\,(U^c) \cap A^* = \varnothing$, which implies that $U^c$ is an attracting neighborhood for $A = \mathrm{Inv}(U^c)$. Thus, by Lemma 3.23, $\mathrm{cl}\,(U) \cap A = \mathrm{cl}\,(U) \cap R^* = \varnothing$.

If $U$ is a neighborhood of $R$ with $\mathrm{cl}\,(U) \cap R^* = \varnothing$, then $\mathrm{cl}\,(U^c) \cap R = \mathrm{cl}\,(U^c) \cap A^* = \varnothing$, and therefore $U^c$ is an attracting neighborhood with $A = \mathrm{Inv}(U^c)$. Consequently, $U$ is a repelling neighborhood for $R$ by Corollary 3.24. ■

COROLLARY 3.26. *If $U \subset X$ is a repelling neighborhood, then also $\mathrm{cl}\,(U)$ and $\mathrm{int}\,(U)$ are repelling neighborhoods and $R = \alpha(U) = \alpha(\mathrm{cl}\,(U)) = \alpha(\mathrm{int}\,(U))$.*

PROOF. The complement $U^c$ is an attracting neighborhood by Corollary 3.24. Lemma 3.2 then yields that also $\mathrm{int}\,(U^c)$ and $\mathrm{cl}\,(U^c)$ are attracting neighborhoods. Consequently, $\mathrm{cl}\,(U) = (\mathrm{int}\,(U^c))^c$ and $\mathrm{int}\,(U) = (\mathrm{cl}\,(U^c))^c$ are repelling neighborhoods. Since $A = R^* = \omega(U^c) = \omega(\mathrm{cl}\,(U^c)) = \omega(\mathrm{int}\,(U^c))$, we have that $\alpha(\mathrm{cl}\,(U)) = A^* = R$ and $\alpha(\mathrm{int}\,(U)) = A^* = R$. ■

COROLLARY 3.27. *Let $U \subset X$ be a repelling neighborhood for a repeller $R = \alpha(U)$, and let $U' \subset X$ be a neighborhood of $R$ such that $R \subset U' \subset \mathrm{cl}\,(U)$. Then also $U'$ is a repelling neighborhood with $R = \alpha(U')$.*

PROOF. Since $\mathrm{cl}\,(U') \subset \mathrm{cl}\,(U)$, we have $\mathrm{cl}\,(U') \cap R^* = \varnothing$. ■

*Proof of Corollary 3.15.* If $U$ is an repelling neighborhood, then $R = \alpha(U) \subset \mathrm{int}\,(U)$ is a repeller, which is compact and invariant. To show that $R$ is the maximal forward invariant set in $U$, suppose that $S$ is a forward invariant set in $U$ satisfying $R \subset S \subset U$. Then, since $R$ is forward-backward invariant,

$$R = \mathrm{cl}\,(R) = \alpha(R) \subset \alpha(S) \subset \alpha(U) = R,$$

and consequently $S \subset \mathrm{cl}\,(S) \subset \alpha(S) = R$ by Proposition 2.13(viii), and thus $R = S$. This proves that $R = \mathrm{Inv}^+(U)$. Using duality, Corollary 3.26 states that $\alpha(U) = \alpha(\mathrm{cl}\,(U))$, and therefore the the same reasoning as above implies that $R = \mathrm{Inv}^+(\mathrm{cl}\,(U))$. ■

PROPOSITION 3.28. *If $R \subset X$ is a repeller and $R' \subset R$ is a repeller for $\varphi|_R$, then $R'$ is a repeller (for $\varphi$).*

PROOF. By definition $R \in \mathsf{Invset}^{\pm}(\varphi)$ and for the dual attractor to $R$ we have $A = R^* \in \mathsf{Invset}(\varphi)$. Let $A'$ be the dual attractor of $R'$ inside $R$ with respect to $\varphi|_R$. Consider the dual set

$$R'^{\ominus} = \{x \in X \mid \exists \gamma_x \text{ such that } \alpha_{\mathrm{o}}(\gamma_x^-) \cap R' = \varnothing\},$$

which, by Proposition 2.16, is invariant. Then, since $R \in \mathsf{Invset}^{\pm}(\varphi)$ and $R' \in \mathsf{Invset}^{\pm}(R)$, it follows that $R' \in \mathsf{Invset}^{\pm}(\varphi)$, and $R'$ is compact. The former follows from the fact that $\varphi|_R(t,x) = \varphi(t,x)$ for all $(t,x) \in \mathbb{T} \times R$.

If $x \in A$, the invariance of $A$ implies that the exists an orbit $\gamma_x \subset A$, and therefore $A \subset R'^{\ominus}$. Similarly, if $x \in A'$, the invariance of $A'$ implies that there exists an orbit $\gamma_x \subset A'$, and consequently $A' \subset R'^{\ominus}$. This implies that

$$A \cup A' \subset R'^{\ominus}.$$

By Proposition 2.16 the forward-backward invariance and compactness of $R'$ implies that $R'^{\ominus} \cap R' = \varnothing$. Therefore we can choose a compact neighborhood $N$ of $R'^{\ominus}$, such that $N \cap R' = \varnothing$. Let $x \in N \setminus R'^{\ominus}$, then $x \in X \setminus R'^{\ominus}$, which implies that $\alpha_{\mathrm{o}}(\gamma_x^-) \cap R' \neq \varnothing$ for all orbits $\gamma_x$, and therefore there are no backward orbits $\gamma_x^- : \mathbb{T}^- \to N$. By Lemma 3.11 we conclude that $R'^{\ominus}$ is an attractor in $X$.

By Proposition 3.16, the dual repeller in $X$ of $R'^{\ominus}$ is $R'^{\ominus *} = \{x \in X \mid \omega(x) \cap R'^{\ominus} = \varnothing\}$. If $x \in R'$, then, since $R'$ is forward-backward invariant in $X$, it holds that $\omega(x) \subset R'$. Since $R' \cap R'^{\ominus} = \varnothing$, we conclude that $R' \subset R'^{\ominus *}$. Suppose $x \in R'^{\ominus *} \subset X \setminus R'^{\ominus}$. If $x \in X \setminus (R'^{\ominus} \cup R)$, then from Theorem 3.19, we have $\omega(x) \subset A$. Similarly, if $x \in R \setminus (R'^{\ominus} \cup R') = R \setminus (A' \cup R')$, then $\omega(x) \subset A'$. In both cases $\omega(x) \cap R'^{\ominus} \neq \varnothing$. On the other hand $R'^{\ominus *}$ is a repeller and therefore is forward-backward invariant, which implies $\omega(x) \subset R'^{\ominus *} \subset X \setminus R'^{\ominus}$, a contradiction. We therefore conclude that $x \in R'$, and thus $R'^{\ominus *} \subset R'$. Combining the inclusions we obtain $R'^{\ominus *} = R'$, which proves that $R'$ is a repeller in $X$. ■

## 4. Lattices of attractors, repellers, and their neighborhoods

Recall that the defining property of an attracting neighborhood $U$ is that $\omega(U) \subset \operatorname{int}(U)$, and for attracting neighborhoods $\omega(U) = \operatorname{Inv}(U)$. The set of attracting neighborhoods $\mathsf{ANbhd}(\varphi)$ has additional structure as a distributive lattice.

PROPOSITION 4.1. *The set* $\mathsf{ANbhd}(\varphi)$ *is a bounded, distributive lattice with respect to the binary operations* $\vee = \cup$ *and* $\wedge = \cap$*. The neutral elements are* $0 = \varnothing$ *and* $1 = X$*.*

PROOF. The elements $\varnothing$ and $X$ are attracting neighborhoods. Let $U, U'$ be attracting neighborhoods for attractors $A, A'$ respectively. Then $U \cup U'$ and $U \cap U'$ are neighborhoods of $A \cup A'$ and $A \cap A'$ respectively. From Proposition 2.11(v) we have

$$\begin{aligned} \omega(U \cup U') &= \omega(U) \cup \omega(U') = A \cup A' \\ &\subset \operatorname{int}(U) \cup \operatorname{int}(U') \subset \operatorname{int}(U \cup U'), \end{aligned} \tag{9}$$

which implies that $U \cup U' \in \mathsf{ANbhd}(\varphi)$.

By Proposition 2.11(iv) and (v), we have

$$\omega(U \cap U') \subset \omega(U) \cap \omega(U') = A \cap A'.$$

From Proposition 2.11(viii), since $\omega(U \cap U')$ is compact and invariant, we deduce that

$$\omega(A \cap A') \subset \omega(U \cap U') = \omega(\omega(U \cap U')) \subset \omega(A \cap A').$$

Hence

$$\omega(U \cap U') = \omega(A \cap A') \subset \operatorname{int}(U) \cap \operatorname{int}(U') = \operatorname{int}(U \cap U'), \tag{10}$$

which proves that $U \cap U' \in \mathsf{ANbhd}(\varphi)$. The distributivity follows since $\mathsf{ANbhd}(\varphi)$ is a sublattice of $\mathsf{P}(X)$ and thus is a lattice of sets. ■

The same arguments can be made for repelling neighborhoods $\mathsf{RNbhd}(\varphi)$ and alpha-limit sets.

PROPOSITION 4.2. *The set* $\mathsf{RNbhd}(\varphi)$ *is a bounded, distributive lattice with respect to the binary operations* $\vee = \cup$ *and* $\wedge = \cap$*. The neutral elements are* $0 = \varnothing$ *and* $1 = X$*.*

PROOF. The elements $\varnothing$ and $X$ are repelling neighborhoods. Let $U, U'$ be repelling neighborhoods for repellers $R, R'$ respectively. Then $U \cup U'$ and $U \cap U'$ are neighborhoods of $R \cup R'$ and $R \cap R'$ respectively. From Proposition 2.13(v) we have

$$\begin{aligned} \alpha(U \cup U') &= \alpha(U) \cup \alpha(U') = R \cup R' \\ &\subset \operatorname{int}(U) \cup \operatorname{int}(U') \subset \operatorname{int}(U \cup U'), \end{aligned} \tag{11}$$

which implies that $U \cup U' \in \mathsf{RNbhd}(\varphi)$.

From Proposition 2.13(viii) we deduce, since $R \cap R' \in \mathsf{Invset}^{\pm}(\varphi)$ and $R \cap R'$ is compact, that $R \cap R' = \alpha(R \cap R') \subset \alpha(U \cap U')$. From Proposition 2.13(iv) and (v) we derive that $\alpha(U \cap U') \subset \alpha(U) \cap \alpha(U') = R \cap R'$ so that $\alpha(U \cap U') = R \cap R'$. Therefore

$$\alpha(U \cap U') = R \cap R' \subset \mathrm{int}\,(U) \cap \mathrm{int}\,(U') = \mathrm{int}\,(U \cap U'), \tag{12}$$

which proves that $U \cap U' \in \mathsf{RNbhd}(\varphi)$. The distributivity follows since $\mathsf{RNbhd}(\varphi)$ is a sublattice of $\mathsf{P}(X)$ and thus is a lattice of sets. ■

Next, using the mapping $\omega : \mathsf{ANbhd}(\varphi) \to \mathsf{Att}(\varphi)$, defined by $U \mapsto \omega(U)$, we establish that $\mathsf{Att}(\varphi)$ is a lattice, and $\omega$ is a lattice homomorphism.

PROPOSITION 4.3. *The set* $\mathsf{Att}(\varphi)$ *is a sublattice of* $\mathsf{Invset}(\varphi)$ *and the mapping* $\omega : \mathsf{ANbhd}(\varphi) \to \mathsf{Att}(\varphi)$ *is surjective homomorphism.*

PROOF. Let $A = \omega(U)$ and $A' = \omega(U')$ be attractors with attracting neighborhoods $U$ and $U'$ respectively. Recall that the lattice operations in $\mathsf{Invset}(\varphi)$ are defined by $\vee = \cup$ and $S \wedge S' = \mathrm{Inv}(S \cap S')$. Since $A, A' \in \mathsf{Invset}^{+}(\varphi)$, we have $A \cap A' \in \mathsf{Invset}^{+}(\varphi)$ so that $\omega(A \cap A') = \mathrm{Inv}(A \cap A')$ by Proposition 2.11(iii). From Equations (9) and (10) we have the identities $A \vee A' = A \cup A' = \omega(U \cup U')$ and $A \wedge A' = \omega(A \cap A') = \mathrm{Inv}(A \cap A') = \omega(U \cap U')$. Therefore $\mathsf{Att}(\varphi)$ is a sublattice of $\mathsf{Invset}(\varphi)$, and $\omega$ is a lattice homomorphism. Also, $\omega(X)$ is the maximal element in $\mathsf{Invset}(\varphi)$ and $\omega(\varnothing) = \varnothing$ is the minimal element in $\mathsf{Invset}(\varphi)$, which proves that $\mathsf{Att}(\varphi)$ is a sublattice and $\omega$ a homomorphism. Surjectivity follows from Proposition 3.5. ■

Note that since $\mathsf{Att}(\varphi)$ is a sublattice of $\mathsf{Invset}(\varphi)$, the attractors are ordered by set inclusion. We have a the same results for $\mathsf{Rep}(\varphi)$ and $\alpha : \mathsf{RNbhd}(\varphi) \to \mathsf{Rep}(\varphi)$.

PROPOSITION 4.4. *The set* $\mathsf{Rep}(\varphi)$ *is a sublattice of* $\mathsf{Invset}^{\pm}(\varphi)$ *and the mapping* $\alpha : \mathsf{RNbhd}(\varphi) \to \mathsf{Rep}(\varphi)$ *is surjective homomorphism.*

PROOF. Let $R = \alpha(U)$ and $R' = \alpha(U')$ be repellers with repelling neighborhoods $U$ and $U'$ respectively. From Equations (11) and (12) we have the identities $R \vee R' = R \cup R' = \alpha(U \cup U')$ and $R \wedge R' = R \cap R' = \alpha(U \cap U')$, which shows that $\mathsf{Rep}(\varphi)$ is a sublattice of $\mathsf{Invset}(\varphi)$. Also, $\alpha(X)$ is the maximal element in $\mathsf{Invset}^{\pm}(\varphi)$ and $\alpha(\varnothing) = \varnothing$ is the minimal element in $\mathsf{Invset}^{\pm}(\varphi)$, which proves that $\mathsf{Rep}(\varphi)$ is a sublattice and $\alpha$ is a homomorphism. Surjectivity follows from Proposition 3.14. ■

REMARK 4.5. Recall now that $\mathrm{Inv}(X) = \omega(X) = 1$ and $\varnothing = 0$ are the neutral elements in $\mathsf{Att}(\varphi)$. For repellers $\mathrm{Inv}^{+}(X) = \alpha(X) = X = 1$ and $\varnothing = 0$ are the neutral elements in $\mathsf{Rep}(\varphi)$. Since $X \in \mathsf{Invset}^{+}(\varphi)$, it follows from Proposition 2.13(viii) that $X \subset \alpha(X) \subset X$ and thus $\alpha(X) = X$.

PROPOSITION 4.6. *The mapping* ${}^{c} : \mathsf{ANbhd}(\varphi) \to \mathsf{RNbhd}(\varphi)$, *defined by* $U \mapsto U^c$ *is a lattice anti-isomorphism.*

PROOF. By Corollary 3.24 the complement of an attracting neighborhood is a repelling neighborhood and vice versa, which implies that the mapping $^c$ is well-defined and bijective. Since the lattice operations are $\cup$ and $\cap$, DeMorgan's laws imply that $^c$ is lattice anti-isomorphism. ■

Define the mappings $^* : \mathsf{Att}(\varphi) \to \mathsf{Rep}(\varphi)$ and $^* : \mathsf{Rep}(\varphi) \to \mathsf{Att}(\varphi)$ by $A = \omega(U) \mapsto A^* = \alpha(U^c)$ and $R = \alpha(U) \mapsto R^* = \omega(U^c)$ respectively.

PROPOSITION 4.7. *The mappings* $^* : \mathsf{Att}(\varphi) \to \mathsf{Rep}(\varphi)$ *and* $^* : \mathsf{Rep}(\varphi) \to \mathsf{Att}(\varphi)$ *are lattice anti-isomorphisms with* $(A^*)^* = A$ *and* $(R^*)^* = R$.

PROOF. From Propositions 4.3 and 4.4

$$\begin{aligned} (A \cup A')^* &= \big(\omega(U) \cup \omega(U')\big)^* = \big(\omega(U \cup U')\big)^* \\ &= \alpha((U \cup U')^c) = \alpha(U^c \cap U'^c) = \alpha(U^c) \cap \alpha(U'^c) \\ &= A^* \cap A'^*. \end{aligned}$$

Similarly

$$\begin{aligned} (A \wedge A')^* &= \big(\omega\big(\omega(U) \cap \omega(U')\big)\big)^* = \big(\omega(U \cap U')\big)^* \\ &= \alpha((U \cap U')^c) = \alpha(U^c \cup U'^c) = \alpha(U^c) \cup \alpha(U'^c) \\ &= A^* \cup A'^*. \end{aligned}$$

Therefore the mapping $^* : \mathsf{Att}(\varphi) \to \mathsf{Rep}(\varphi)$ is a lattice anti-homomorphism. Moreover, it is an anti-isomorphism by Proposition 3.18. The proof for $^* : \mathsf{Rep}(\varphi) \to \mathsf{Att}(\varphi)$ is analogous. ■

The above propositions are summarized in the commutative diagram (1).

## 5. Lifting attractor and repeller lattices

DEFINITION 5.1. Let $\mathsf{K}, \mathsf{L}$ be bounded, distributive lattices, let $h\colon \mathsf{K} \to \mathsf{L}$ be a lattice homomorphism and let $\ell\colon \mathsf{L}' \to \mathsf{L}$ be a lattice homomorphism. A lattice homomorphism $k\colon \mathsf{L}' \to \mathsf{K}$ is a *lift* of $\ell$ through $h$, if $h \circ k = \ell$ or equivalently the following diagram commutes:

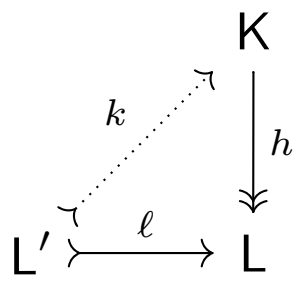

(13)

Observe that using this language the result of Theorem 1.2 can be restated as the existence of the lift of the embedding of a finite sublattice in $\mathsf{Att}(\varphi)$ or $\mathsf{Rep}(\varphi)$ through Inv or Inv$^+$.

Using Birkhoff's Representation Theorem we recast the lifting diagram (13) as

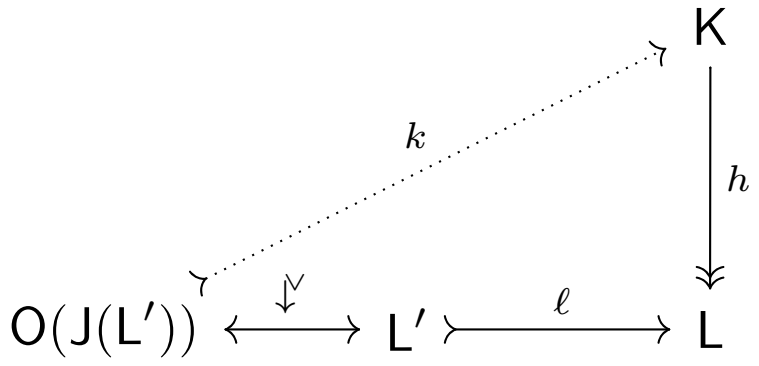

(14)

Observe that since $\downarrow^\vee$ is an isomorphism, if we prove the existence of $k$ which makes the (14) commute, then $k \circ \downarrow^\vee$ provides a lift of $\ell$. To simplify the notation, in what follows let $\mathsf{P} = \mathsf{J}(\mathsf{L}')$ and let $s = \ell \circ (\downarrow^\vee)^{-1}$. In this language our goal is to find a lift $k$ of $s$ through $h$:

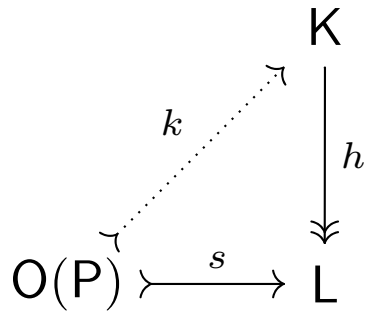

(15)

REMARK 5.2. Let $\mathsf{P}$ be a poset. Let $\lambda \in \mathsf{O}(\mathsf{P})$. By definition $\lambda$ is a down set in $\mathsf{P}$ and hence can simultaneously be viewed as a subset of $\mathsf{P}$.

LEMMA 5.3. *Let $a, b, c$ be elements in a Boolean algebra $\mathsf{B}$. Then*

$$c \wedge b^c \wedge a = 0 \qquad \Longleftrightarrow \qquad c \wedge a \leq b. \tag{16}$$

PROOF. $0 = c \wedge b^c \wedge a = (c \wedge a) \wedge b^c$, which is equivalent to $c \wedge a \leq b$ by (2). ■

Observe that if $k\colon \mathsf{O}(\mathsf{P}) \to \mathsf{K}$ is a lift, and $\alpha, \beta, \gamma \in \mathsf{O}(\mathsf{P})$, then (16) can be rewritten as

$$\gamma \cap \alpha \subset \beta \quad \Leftrightarrow \quad k(\gamma) \wedge k(\alpha) \leq k(\beta). \tag{17}$$

Let $\lambda \in \mathsf{O}(\mathsf{P})$. Note that $0 \in \mathsf{O}(\lambda)$. However, if $\lambda \neq \mathsf{P}$, then $\mathsf{P} \notin \mathsf{O}(\lambda)$, and hence $\mathsf{O}(\lambda)$ is not a sublattice of $\mathsf{O}(\mathsf{P})$. Therefore we define $\lambda^\top$ to be the poset $\lambda \cup \{\top\}$ where the additional top element $\top$ has relations $p \leq \top$ for all $p \in \lambda$. Then

$$\mathsf{O}(\lambda^\top) \approx \{\alpha \in \mathsf{O}(\mathsf{P}) \mid \alpha \subset \lambda \text{ or } \alpha = \mathsf{P}\}$$

making $\mathsf{O}(\lambda^\top)$ a sublattice of $\mathsf{O}(\mathsf{P})$. This implies that the Booleanization $\mathsf{B}(\mathsf{O}(\lambda^\top)) \subset \mathsf{B}(\mathsf{O}(\mathsf{P})) = 2^{\mathsf{P}}$.

DEFINITION 5.4. Let $\lambda \in \mathsf{O}(\mathsf{P})$. A lattice homomorphism $k : \mathsf{O}(\lambda^\top) \to \mathsf{K}$ is a *partial lift of $s$ on* $\mathsf{O}(\lambda^\top)$ if

$$h(k(\beta)) = s(\beta) \text{ for all } \beta \leq \lambda.$$

Note that by the above definition $k(1) = 1$, since $k$ is a lattice homomorphism.

DEFINITION 5.5. A partial lift is a *conditional lift of $s$ on $\mathsf{O}(\lambda^\top)$* if there exists $\{v_\alpha \in h^{-1}(s(\alpha)) \mid \alpha \in \mathsf{O}(\mathsf{P})\}$ such that

$$k(\gamma) \wedge v_\alpha \leq k(\beta) \tag{18}$$

for all $\beta, \gamma \subset \lambda$ for which $\gamma \cap \alpha \subset \beta$. The elements $v_\alpha$ are called *conditioners* for the partial lift $k : \mathsf{O}(\lambda^\top) \to \mathsf{K}$.

REMARK 5.6. Let $k$ be a conditional lift of $s$ on $\mathsf{O}(\lambda^\top)$. The following statements indicate possible choices for conditioners.

(1) If $\alpha \subset \lambda$, then $k(\alpha)$ is a conditioner.
(2) If $\gamma \cap \alpha \not\subset \beta$, then any $v_\alpha \in h^{-1}(s(\alpha))$ is a conditioner.
(3) Assume $v_\alpha \in h^{-1}(s(\alpha))$ is a conditioner. If $v'_\alpha \in h^{-1}(s(\alpha))$ and $v'_\alpha \leq v_\alpha$, then $v'_\alpha$ is a conditioner.

We make use of the following equivalent characterization of a conditional lift.

PROPOSITION 5.7. *Let $k\colon \mathsf{O}(\lambda^\top) \to \mathsf{K}$ be a partial lift. Then, $k$ is a conditional lift with conditioners $\{v_\alpha \in h^{-1}(s(\alpha)) \mid \alpha \in \mathsf{O}(\mathsf{P})\}$ if and only if for any $p \in \lambda$ and any $\alpha \in \mathsf{O}(\mathsf{P})$ such that $p \notin \alpha$*

$$\mathsf{B}(k)(\{p\}) \wedge v_\alpha = 0.$$

PROOF. If $k$ is a conditional lift, then we have the existence of a set $\{v_\alpha \in h^{-1}(s(\alpha)) \mid \alpha \in \mathsf{O}(\mathsf{P})\}$ on which (18) is satisfied. Let $\gamma = \downarrow p$ and let $\beta = \overleftarrow{\gamma}$. Observe that $\beta, \gamma \subset \lambda$ and $\gamma \cap \alpha \subset \beta$, thus

$$\begin{aligned} k(\gamma) \wedge v_\alpha \leq k(\beta) &\Leftrightarrow k(\gamma) \wedge k(\beta)^c \wedge v_\alpha = 0 \quad \text{by Lemma 5.3} \\ &\Leftrightarrow \mathsf{B}(k)(\gamma \wedge \beta^c) \wedge v_\alpha = 0 \\ &\Rightarrow \mathsf{B}(k)(\{p\}) \wedge v_\alpha = 0. \end{aligned} \tag{19}$$

To prove the converse, consider $\beta, \gamma \subset \lambda$ such that $\gamma \cap \alpha \subset \beta$. Observe that

$$\begin{aligned} \mathsf{B}(k)(\gamma \cap \beta^c) &= \bigvee_{p \in \gamma \cap \beta^c} \mathsf{B}(k)(\{p\}) \\ \mathsf{B}(k)(\gamma \cap \beta^c) \wedge v_\alpha &= \bigvee_{p \in \gamma \cap \beta^c} \mathsf{B}(k)(\{p\}) \wedge v_\alpha = 0 \end{aligned}$$

and we can now return to the equivalences of (19). ■

DEFINITION 5.8. Let $\mathsf{K}, \mathsf{L}$ be bounded, distributive lattices. A lattice epimorphism $h : \mathsf{K} \to \mathsf{L}$ is *spacious*, if the following two conditions hold for every lattice embedding $s : \mathsf{O}(\mathsf{P}) \to \mathsf{L}$.

(i) For every minimal $q \in \mathsf{P}$, every partial lift $k : \mathsf{O}(\{q\}^\top) \to \mathsf{K}$ of $s$ on $\mathsf{O}(\{q\}^\top)$ is a conditional lift on $\mathsf{O}(\{q\}^\top)$.

(ii) For every $\lambda \in \mathsf{O}(\mathsf{P})$ and every partial lift $k\colon \mathsf{O}(\lambda^{\top}) \to \mathsf{K}$ of $s$ on $\mathsf{O}(\lambda^{\top})$ and any minimal $q \in \mathsf{P} \setminus \lambda$, there exists

$$\big\{v_\alpha \in h^{-1}(s(\alpha)) \mid \alpha \in \mathsf{O}(\mathsf{P}) \text{ such that } q \notin \alpha\big\}$$

for which

$$v_\mu \wedge v_\alpha \leq k(\lambda) \quad \text{where } \mu :=\downarrow q. \tag{20}$$

PROPOSITION 5.9. *Let $h : \mathsf{K} \to \mathsf{L}$ be a surjective, lattice homomorphism between bounded, distributive lattices. If $h^{-1}(0) = 0$, then Condition (i) of Definition 5.8 is satisfied.*

PROOF. For any $\alpha \in \mathsf{P}$, choose $v_\alpha \in h^{-1}(s(\alpha))$. We remark that this implies that $v_0 = 0$. Let $s\colon \mathsf{O}(\mathsf{P}) \to \mathsf{L}$ be a lattice embedding, and let $q \in \mathsf{P}$ be minimal. Let $k\colon \mathsf{O}(\{q\}) \to \mathsf{K}$ be a partial lift of $s$. Since $q$ is minimal, if $\beta, \gamma \subset \{q\}$, then $\gamma, \beta \in \{\varnothing, \{q\}\}$. If $\gamma = \{q\}$ and $\beta = \varnothing$, then

$$h(k(\{q\}) \wedge v_\alpha) = h(k(\{q\})) \wedge h(v_\alpha) = s(\{q\}) \wedge s(\alpha) = s\left(\{q\} \cap \alpha)\right) = s(0) = 0$$

for all $q \notin \alpha$. Since $h^{-1}(0) = 0$, this implies that $k(\{q\}) \wedge v_\alpha = 0$, for all $q \notin \alpha$. If $\gamma = \beta = \{q\}$, then

$$k(\gamma) \wedge v_\alpha = k(\{q\}) \wedge v_\alpha \leq k(\{q\}) = k(\beta),$$

and if $\gamma = \beta = \varnothing$, then

$$k(\gamma) \wedge v_\alpha = k(0) \wedge v_\alpha = 0 \wedge v_\alpha = 0 = k(\beta),$$

which proves the lemma. ■

COROLLARY 5.10. $\mathrm{Inv}\colon \mathsf{ANbhd}(\varphi) \to \mathsf{Att}(\varphi)$ *satisfies Condition (i) of Definition 5.8. If $\varphi$ is surjective, then* $\mathrm{Inv}^+\colon \mathsf{RNbhd}(\varphi) \to \mathsf{Rep}(\varphi)$ *satisfies Condition (i) of Definition 5.8.*

PROOF. By Proposition 2.11(ii) if $U \in \mathsf{ANbhd}(X)$ and $U \neq \varnothing$, then $\mathrm{Inv}(U, \varphi) \neq \varnothing$. By Proposition 2.13(ii) if $\varphi$ is surjective, then $U \in \mathsf{RNbhd}(X)$ and $U \neq \varnothing$ implies $\mathrm{Inv}^+(U, \varphi) \neq \varnothing$. ■

EXAMPLE 5.11. To show that $\mathrm{Inv}\colon \mathsf{ANbhd}(X) \to \mathsf{Att}(X)$ is not in general spacious, we consider a semi-flow $\varphi$ on a 1-dimensional graph as shown in Figure 1[left] with fixed points $\{1, 2, 3\}$. The lattice of attractors $\mathsf{Att}(\varphi)$ is indicated in Figure 1[right]. Consider $\mathsf{L}' = \mathsf{Att}(\varphi)$. Then $\mathsf{J}(\mathsf{L}') = \mathsf{P} = \{1, 2, 3\}$ with the order described in Figure 1[middle].

Choose $\lambda = \{1\}$. Observe that $q = 2$ is a minimal element in $\mathsf{P} \setminus \lambda = \{2, 3\}$. Let $\mu =\downarrow q =\downarrow 2 = \{1, 2\}$. Consider the down set $\alpha = \{1, 3\}$, and note that $2 \notin \alpha$. Furthermore,

$$\mu \wedge \alpha = \{1, 2\} \cap \{1, 3\} = \{1\}.$$

Let $k\colon \mathsf{O}(\{1\}) \to \mathsf{ANbhd}(X)$ be a partial lift such that $k(\{1\})$ is the attracting neighborhood for $s(\{1\})$ indicated in Figure 1 by a thick line. For every choice of attracting neighborhoods $V_{\{1,2\}}$ and $V_{\{1,3\}}$ of $s(\{1,2\})$ and $s(\{1,3\})$ respectively, we have that

$$k(\{1\}) \subsetneq V_{\{1,2\}} \cap V_{\{1,3\}} = s(\{1,2\}) \wedge s(\{1,3\}),$$

contradicting (20). Since the spacious condition has to be satisfied for all partial lifts, we conclude that $\mathrm{Inv} : \mathsf{ANbhd}(X) \to \mathsf{Att}(X)$ is not spacious in general.

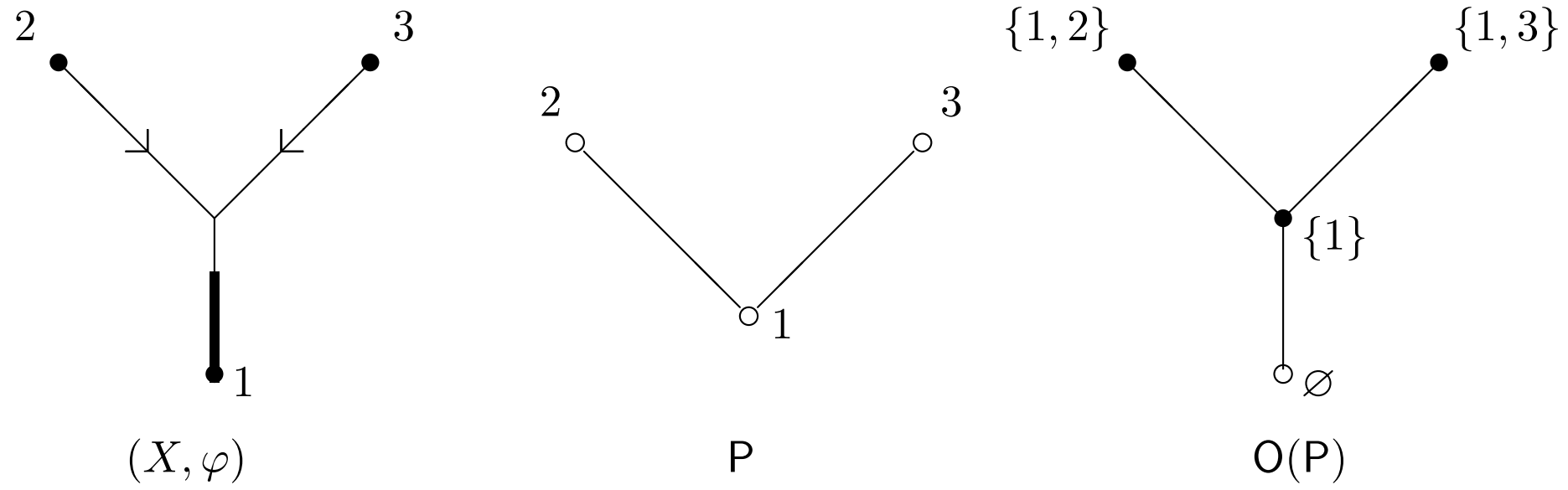

FIGURE 1. A 1-dimensional semi-flow [left], the order of fixed points [middle], and the lattice of attractors [right].

PROPOSITION 5.12. $\mathrm{Inv}^+\colon \mathsf{RNbhd}(\varphi) \to \mathsf{Rep}(\varphi)$ *is spacious.*

PROOF. Let $s\colon \mathsf{O}(\mathsf{P}) \to \mathsf{Rep}(X)$ be a lattice embedding for some finite poset $\mathsf{P}$. Let $q \in \mathsf{P}$ be minimal. Assume $k\colon \mathsf{O}(\{q\}^{\top}) \to \mathsf{RNbhd}(\varphi)$ is a partial lift of $s$. We need to show that $k$ is a conditional lift on $\mathsf{O}(\{q\}^{\top})$.

We make two observations concerning $\alpha \in \mathsf{O}(\mathsf{P})$ for which $q \notin \alpha$. First, via (17)

$$s(\{q\}) \wedge s(\alpha) = s(\{q\} \cap \alpha) = s(0) = 0,$$

and second, there exists a compact neighborhood $N_\alpha \subset X$ of $s(\alpha)$ such that $k(\{q\}) \cap N_\alpha = \varnothing$. Let $W_\alpha \in \mathsf{RNbhd}(\varphi)$ be a repelling neighborhood for $s(\alpha)$. By Corollary 3.27, $v_\alpha := W_\alpha \cap N_\alpha$ is a repelling neighborhood of $s(\alpha)$, i.e. $\mathrm{Inv}^+(v_\alpha) = s(\alpha)$, and furthermore

$$k(\{q\}) \wedge v_\alpha = 0.$$

Since $q$ is minimal, if $\beta, \gamma \subset \{q\}$, then $\gamma, \beta \in \{\varnothing, \{q\}\}$. If $\gamma = \{q\}$ and $\beta = \varnothing$, then

$$k(\{q\}) \wedge v_\alpha = 0 = k(0) = k(\beta).$$

for all $q \notin \alpha$. If $\gamma = \beta = \{q\}$, then

$$k(\gamma) \wedge v_\alpha = k(\{q\}) \wedge v_\alpha = 0 \leq k(\beta).$$

If $\gamma = \beta = \varnothing$, then

$$k(\gamma) \wedge v_\alpha = k(0) \wedge v_\alpha = 0 \wedge v_\alpha = 0 = k(\beta).$$

Thus the partial lift $k$ satisfies Definition 5.8(i). We now verify Definition 5.8(ii). Let $\lambda \in \mathsf{O}(\mathsf{P})$ and let $q \in \mathsf{P} \setminus \lambda$ be minimal. Define $\mu = \downarrow q$. Let $k\colon \mathsf{O}(\lambda^\top) \to \mathsf{RNbhd}(\varphi)$ be a partial lift.

As above let $\alpha \in \mathsf{O}(\mathsf{P})$ for which $q \notin \alpha$. Observe that $\mu \cap \lambda^c \cap \alpha = \varnothing$, hence by (17) $\mu \cap \alpha \subset \lambda$. In particular,

$$s(\mu \cap \alpha) \leq s(\lambda) \subset \text{int}\,(k(\lambda)).$$

Since $[\text{int}\,(k(\lambda))]^c$ and $s(\mu \cap \alpha)$ are disjoint compact sets there exists $\epsilon(\alpha) > 0$ such that

$$\text{dist}\,(s(\mu \cap \alpha), [\text{int}\,(k(\lambda))]^c) = \epsilon(\alpha).$$

Define $\epsilon := \min\{\epsilon(\alpha) \mid q \notin \alpha\}$. Choose $\{W_\xi \mid \xi \in \mathsf{O}(\mathsf{P})\} \subset \mathsf{RNbhd}(\varphi)$ such that $W_\xi$ is a repelling neighborhood for $s(\xi)$. By Remark 5.6 we can define the conditioners as follows:

$$v_\xi := \begin{cases} k(\xi) & \text{if } \xi \leq \lambda, \\ W_\xi \cap B_{\epsilon/2}(s(\xi)) & \text{if } q \notin \xi \text{ or } \xi = \mu, \\ W_\xi & \text{if } q \in \xi \text{ and } \xi \neq \mu. \end{cases}$$

By Corollary 3.27, $v_\xi \in \mathsf{RNbhd}(\varphi)$ for all $\xi \in \mathsf{O}(\mathsf{P})$. By construction, $v_\mu \cap v_\alpha \subset \text{int}\,(k(\lambda)) \subset k(\lambda)$, and hence

$$v_\mu \wedge v_\alpha \leq k(\lambda).$$

Thus (20) is satisfied and hence $\text{Inv}^+\colon \mathsf{RNbhd}(X) \to \mathsf{Rep}(X)$ is spacious. ∎

For spacious lattice homomorphisms we obtain the following lifting theorem.

THEOREM 5.13. *Let* $\mathsf{K}$ *and* $\mathsf{L}$ *be a bounded, distributive lattices, and let* $h\colon \mathsf{K} \to \mathsf{L}$ *be a lattice epimomorphism. If h is spacious and* $h^{-1}(1) = 1$*, then every lattice embedding* $s\colon \mathsf{O}(\mathsf{P}) \to \mathsf{L}$*, with* $\mathsf{P}$ *is finite, admits a lift.*

PROOF. We construct the desired lift $k\colon \mathsf{O}(\mathsf{P}) \to \mathsf{K}$ inductively. To begin set $k(0) = 0$. Let $q \in \mathsf{P}$ be minimal. Define $\lambda = \downarrow q = \{q\} \in \mathsf{O}(\mathsf{P})$. Choose $k(\lambda) \in h^{-1}(s(\lambda))$. Observe that $k\colon \mathsf{O}(\lambda^\top) \to \mathsf{K}$ is a partial lift. Condition (i) of Definition 5.8 implies that $k$ is conditional.

To perform the inductive step, let $\lambda \in \mathsf{O}(\mathsf{P})$ and assume $k\colon \mathsf{O}(\lambda^\top) \to \mathsf{K}$ is a conditional lift of $s$. Let $\left\{v^0_\alpha \in h^{-1}(s(\alpha)) \mid \alpha \in \mathsf{O}(\mathsf{P})\right\}$ denote the conditioners of $k$. Choose a minimal element $q \in \mathsf{P} \setminus \lambda$ and let $\mu = \downarrow q \in \mathsf{O}(\mathsf{P})$. The goal is to extend $k$ to a partial lift $k\colon \mathsf{O}((\lambda \cup \mu)^\top) \to \mathsf{K}$ that is conditional.

Applying the Booleanization functor, we obtain $\mathsf{B}(k)\colon \mathsf{B}(\mathsf{O}(\lambda^\top)) \to \mathsf{B}(\mathsf{K})$. By Proposition 2.3 we have

$$k(\alpha) = \begin{cases} \bigvee_{p\in\alpha} \mathsf{B}(k)(\{p\}) & \text{for all } \alpha \in \mathsf{O}(\lambda^\top), \alpha \neq \mathsf{P} \\ 1 & \text{for } \alpha = \mathsf{P}. \end{cases}$$

We want to extend to a map on $\mathsf{O}((\lambda \cup \mu)^\top)$ by determining a proper image for $\{q\}$ in $\mathsf{B}(\mathsf{K})$. Since $h$ is spacious, Definition 5.8(ii) implies that there exists $\left\{v^1_\alpha \in h^{-1}(s(\alpha)) \mid \alpha \in \mathsf{O}(\mathsf{P})\right\}$ such that $v^1_\mu \wedge v^1_\alpha \leq k(\lambda)$ whenever $q \notin \alpha$. Define

$$v_\alpha = v^0_\alpha \wedge v^1_\alpha \ \text{ for all } \alpha \in \mathsf{O}(\mathsf{P}).$$

Now define

$$B_q = v_\mu \wedge k(\lambda)^c \in \mathsf{B}(\mathsf{K}).$$

Finally define

$$k(\alpha) = \begin{cases} \bigvee_{p\in\alpha} \mathsf{B}(k)(\{p\}) & \text{for all } \alpha \in \mathsf{O}((\lambda\cup\mu)^\top) \text{ with } q \notin \alpha, \alpha \neq \mathsf{P} \\ B_q \vee \bigvee_{p\in\alpha, p\neq q} \mathsf{B}(k)(\{p\}) & \text{for all } \alpha \in \mathsf{O}((\lambda\cup\mu)^\top) \text{ with } q \in \alpha, \alpha \neq \mathsf{P} \\ 1 & \text{for } \alpha = \mathsf{P}. \end{cases}$$

This is a well-defined mapping that extends the domain of $k : \mathsf{O}(\lambda^\top) \to \mathsf{B}(\mathsf{K})$ to $\mathsf{O}((\lambda \cup \mu)^\top)$. The proof is complete once it is shown that $k$ is a conditional partial lift of $s$ on $\mathsf{O}((\lambda \cup \mu)^\top)$. There are three properties of $k$ that need to be shown: (a) $k$ is a lattice homomorphism, (b) $k$ maps into $\mathsf{K}$ with $h \circ k = s$, and (c) $\{v_\alpha\}$ are conditioners for $k$.

*Proof of* (a): For notational convenience, let $B_p = \mathsf{B}(k)(\{p\})$ for $p \in \lambda$ so that

$$k(\alpha) = \bigvee_{p\in\alpha} B_p \ \text{ for all } \alpha \in \mathsf{O}((\lambda\cup\mu)^\top). \tag{21}$$

Since $v_\mu \wedge v_\alpha \leq k_\lambda$ is equivalent to $v_\mu \wedge k(\lambda)^c \wedge v_\alpha = 0$, we have $B_q \wedge v_\alpha = 0$ for all $\alpha \not\ni q$ by definition of $B_q$. Therefore

$$B_p \wedge v_\alpha = 0 \text{ for all } p \in (\lambda\cup\mu) \setminus \alpha. \tag{22}$$

Since $B_q \wedge k(\lambda) = 0$, it follows from equation (4) in Proposition 2.3 that

$$B_p \wedge B_{p'} = 0 \ \text{ for all distinct elements } p, p' \in \lambda \cup \mu. \tag{23}$$

From (21) we have that

$$k(\alpha) \vee k(\alpha') = \Big(\bigvee_{p\in\alpha} B_p\Big) \vee \Big(\bigvee_{p'\in\alpha'} B_{p'}\Big) = \Big(\bigvee_{p\in\alpha\cup\alpha'} B_p\Big) = k(\alpha\cup\alpha').$$

Moreover from (23)

$$\begin{aligned} k(\alpha) \wedge k(\alpha') &= \Big(\bigvee_{p\in\alpha} B_p\Big) \wedge \Big(\bigvee_{p'\in\alpha'} B_{p'}\Big) \\ &= \bigvee_{(p,p')\in\alpha\times\alpha'} (B_p \wedge B_{p'}) = \bigvee_{p\in\alpha\cap\alpha'} B_p = k(\alpha\cap\alpha'), \end{aligned}$$

which proves that $k : \mathsf{O}((\lambda\cup\mu)^{\top}) \to \mathsf{B}(\mathsf{K})$ is a lattice homomorphism.

*Proof of* (b): Setting $\alpha = \mu$ in Equation (22), $B_p \wedge v_\mu = 0$ for all $p \in \lambda \setminus \mu$. Since $\mu =\downarrow q$, it has an immediate predecessor $\overleftarrow{\mu} = \mu \setminus \{q\}$, and since $q \notin \lambda$, we have $\lambda \wedge \mu^c = \lambda \wedge (\overleftarrow{\mu})^c$.

Now

$$k(\lambda) \wedge k(\overleftarrow{\mu})^c = k(\lambda \wedge (\overleftarrow{\mu})^c) = \bigvee_{p\in\lambda\cap\overleftarrow{\mu}^c} B_p = \bigvee_{p\in\lambda\cap\mu^c} B_p.$$

Therefore $k(\lambda) \wedge k(\overleftarrow{\mu})^c \wedge v_\mu = 0$ so that $k(\overleftarrow{\mu}) \vee (k(\lambda) \wedge v_\mu) = k(\overleftarrow{\mu})$. Hence

$$\begin{aligned} k(\mu) &= k(\overleftarrow{\mu}) \vee B_q = k(\overleftarrow{\mu}) \vee (v_\mu \wedge k(\lambda)^c) \\ &= k(\overleftarrow{\mu}) \vee (v_\mu \wedge k(\lambda)) \vee (v_\mu \wedge k(\lambda)^c) \\ &= k(\overleftarrow{\mu}) \vee v_\mu. \end{aligned}$$

Note that $\overleftarrow{\mu} \in \mathsf{O}(\lambda^{\top})$ since $q$ is a minimal element in $\mathsf{P}\setminus\lambda$. This implies that $k(\mu) \in \mathsf{K}$ because both $k(\overleftarrow{\mu})$ and $v_\mu$ are in $\mathsf{K}$. Finally $h(k(\mu)) = h(k(\overleftarrow{\mu})) \vee h(v_\mu) = s(\overleftarrow{\mu}) \vee s(\mu) = s(\mu)$.

Since $\mu =\downarrow q$, the elements $\alpha \in \mathsf{O}((\lambda\cup\mu)^{\top})$ satisfy either $\alpha\cap\mu = 0$, which implies $\alpha \subset \lambda$, or $\alpha = \alpha' \cup \mu$ for some $\alpha' \subset \lambda$. In the former case, we already have $k(\alpha) \in \mathsf{K}$ and $h(k(\alpha)) = s(\alpha)$, since $k$ is a partial lift on $\mathsf{O}(\lambda^{\top})$. Now consider $\alpha = \alpha' \cup \mu$ with $\alpha' \subset \lambda$. Since $k$ is a homomorphism, $k(\alpha) = k(\alpha') \vee k(\mu)$ so that $k(\alpha) \in \mathsf{K}$ and $h(k(\alpha)) = h(k(\alpha')) \vee h(k(\mu)) = s(\alpha') \vee s(\mu) = s(\alpha)$. We have now established $h \circ k = s$ on $\mathsf{O}((\lambda\cup\mu)^{\top})$. There for $k : \mathsf{O}((\lambda\cup\mu)^{\top}) \to \mathsf{K}$ is a partial lift of $s$.

Now that we have established that $k$ on $\mathsf{O}((\lambda\cup\mu)^{\top})$ is an extension of the original lattice homomorphism on $\mathsf{O}(\lambda^{\top})$, and that both domains are sublattices of $\mathsf{O}(\mathsf{P})$, then the functoriality of the Booleanization $\mathsf{B}$ guarantees that $\mathsf{B}(k)(\{p\})$ has not changed for $p \in \lambda$ and $\mathsf{B}(k)(\{q\}) = B_q$, which all lie in $\mathsf{B}(\mathsf{K})$.

*Proof of* (c): This partial lift is conditional with respect to $\{v_\alpha\}$ by Equation (22) using the equivalence in (19). We have now shown that any conditional lift on $\mathsf{O}(\lambda^{\top})$ can be extended to a conditional lift on $\mathsf{O}((\lambda\cup\mu)^{\top})$, where $\mu =\downarrow q$ for some minimal element $q$ of $\mathsf{P} \setminus \lambda$.

We construct consecutive conditional lifts by depleting the set $\mathsf{P}$ as described in the above induction step. This procedure terminates after finitely many steps. To complete the proof we need to show that $k(1) = 1$ complies with the terminal step. Indeed at the terminal step $\lambda \cup \mu = \mathsf{P}$, and hence

$$h(k(1)) = h(k(\lambda)) \vee h(k(\mu)) = s(\lambda) \vee s(\mu) = s(\lambda \cup \mu) = s(1) = 1,$$

and therefore $k(1) = 1$ since $h^{-1}(1) = 1$. ∎

Note that if the hypothesis $h^{-1}(1) = 1$ is omitted, then the above proof provides a lift mapping that satisfies all properties of a lattice homomorphism except possibly $k(1) = 1$.

*Proof of Theorem 1.2.* The proof is done in two steps. We begin by verifying (ii). Let $\mathsf{R}$ be a finite sublattice of $\mathsf{Rep}(\varphi)$. As is indicated in the discussion preceding (15), it is sufficient to lift $s_R\colon \mathsf{O}(\mathsf{P}) \to \mathsf{Rep}(\varphi)$ where $s_R = i \circ (\downarrow^\vee)^{-1}$ and $i$ is the inclusion map. By Proposition 5.12, $\mathrm{Inv}^+$ is spacious. Also, since $X$ is a repeller, $(\mathrm{Inv}^+)^{-1}(X) = X$, see Remark 4.5. Hence by Theorem 5.13, the desired lift $k_R\colon \mathsf{O}(\mathsf{P}) \to \mathsf{RNbhd}(\varphi)$ exists.

As in the proof of (ii), we begin the proof of (i) by considering a finite sublattice $\mathsf{A}$ of $\mathsf{Att}(\varphi)$, and we note that it is sufficient to lift $s_A\colon \mathsf{O}(\mathsf{P}) \to \mathsf{Att}(\varphi)$ where $s_A = i \circ (\downarrow^\vee)^{-1}$ and $i$ is the inclusion map. The remainder of the proof is based on the following diagram that we demonstrate is commutative.

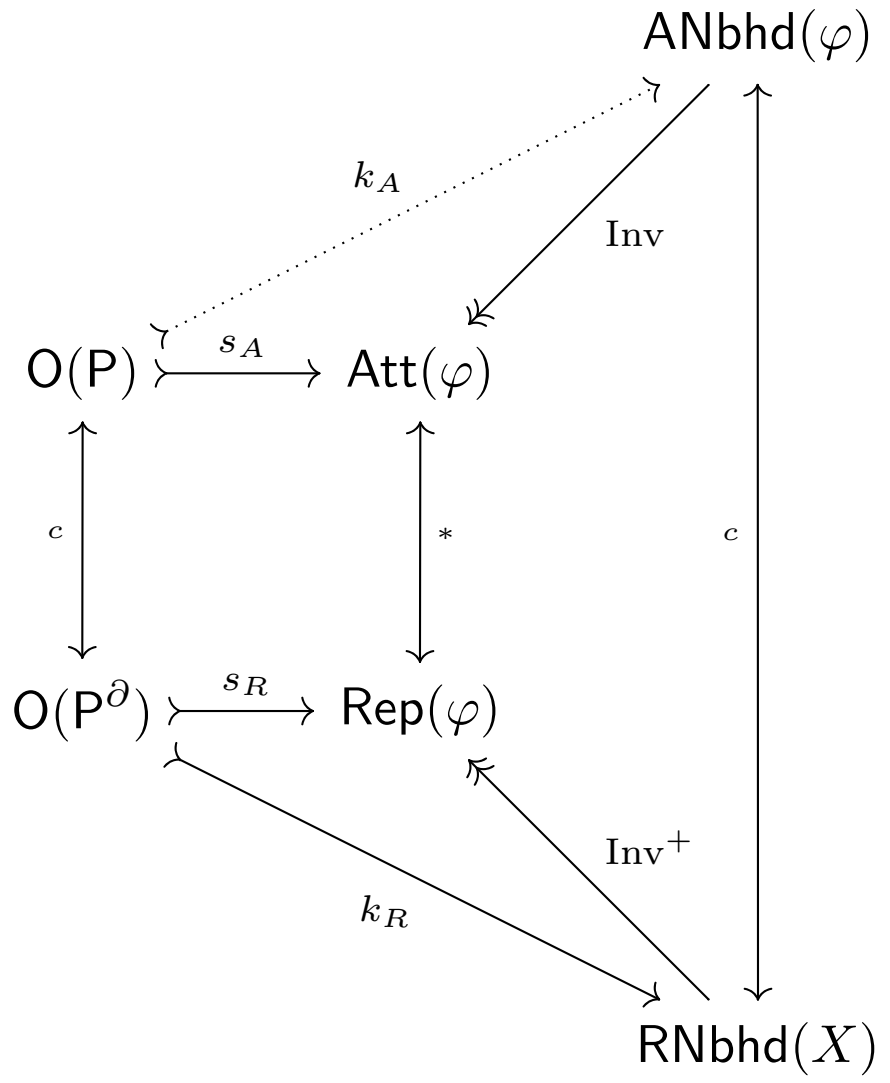

(24)

The complement map ${}^c\colon \mathsf{O}(\mathsf{P}) \to \mathsf{O}(\mathsf{P}^\partial)$ is given by (3). The vertical double headed arrows are involutions. The right square is given by (1) and hence is commutative. Define $s_R\colon \mathsf{O}(\mathsf{P}^\partial) \to \mathsf{Rep}(\varphi)$ by

$${}^* \circ s_A \circ {}^c$$

so that the left square commutes. By Theorem 1.2(ii), there exists $k_R$ such that the lower triangle commutes. Define

$$k_A := {}^c \circ k_R \circ {}^c.$$

The commutativity of the diagram guarantees that $\mathrm{Inv} \circ k_A = s_A$ and hence that $k_A$ is a lift of $s_A$. ■

mischaik@math.rutgers.edu
wkalies@fau.edu
vdvorst@few.vu.nl